\documentclass[11pt,a4paper]{amsart}
\usepackage{a4wide}
\usepackage{amsfonts,amsthm,amsmath,amssymb}
\usepackage{amscd,color}
\usepackage{verbatim}
\usepackage{graphicx,subfigure}
\usepackage{pgfplots} 
\usepackage{enumitem}
\usepackage[normalem]{ulem} 

\newtheorem*{thmA}{Theorem A}
\newtheorem*{thmB}{Theorem B}

\theoremstyle{plain}

\newtheorem{theorem}{Theorem}[section]
\newtheorem{lemma}[theorem]{Lemma}
\newtheorem{corollary}[theorem]{Corollary}
\newtheorem{prop}[theorem]{Proposition}
\newtheorem{proposition}[theorem]{Proposition}

\newcommand{\Hfam}{H_{d}}

\newcommand{\Tfam}{T_{d}}

\newsavebox{\savepar}

\usepackage[normalem]{ulem}

\begin{document}

\title[Boundedness \& simple connectivity of basins of attraction]{Boundedness and simple connectivity of the basins of attraction for some numerical methods}

\date{}

\author{Jordi Canela}
\address{Institut Universitari de Matem\`atiques i Aplicacions de Castell\'o and Departament de Matem\`atiques, Universitat Jaume I, 12071 Castell\'o de la Plana, Spain}
\email{canela@uji.es}
\author{Antonio Garijo}
\address{Departament d'Enginyeria Inform\`atica i Matem\`atiques,
Universitat Rovira i Virgili, 43007 Tarragona, Catalonia, Spain}
\email{antonio.garijo@urv.cat}
\author{Xavier Jarque}
\address{Departament de Matem\`atiques i Inform\`atica, Universitat de Barcelona, Gran Via, 585, 08007 Barcelona, Catalonia and  Centre de Recerca Matemàtica, Edifici C, Campus Bellaterra, 08193 Bellaterra, Catalonia}
\email{xavier.jarque@ub.edu}

\thanks{The first and third authors are supported by the project CIGE/2023/102 funded by Generalitat Valenciana and by the project PID2023-147252NB-I00 financed by MICIU/AEI MCIN/AEI/10.13039/501100011033 and by FEDER, UE. The first author is also supported by the project UJI-B2022-46 funded by Universitat Jaume I. The second author is supported by the project  2021SGR-633 funded by Generalitat de Catalunya. The third author is also supported by the Spanish State Research Agency, through the Severo Ochoa and María de Maeztu Program for Centers and Units of Excellence in R\&D (CEX2020-001084-M)}

\begin{abstract}
In this paper we study the dynamics of Halley's and Traub's  root-finding  algorithms applied to a symmetric family of polynomials of degree $d+1\geq 3$.  We discuss the (un)boundedness and simple connectivity of the immediate basins of attraction of the fixed points associated to the roots of the polynomials.  In particular, we show the existence of polynomials for which the immediate basin of attraction of a root is bounded under Halley's method.\newline

\end{abstract}

\maketitle

\section{Introduction}

The study of discrete dynamical systems generated by the iterates of root-finding algorithms has been a cornerstone for understanding these algorithms not only as useful local tools for finding the desired solutions of non-linear equations, but also as a means to determine good initial conditions for efficiently finding all solutions. It was Cayley \cite{Cayley} who first considered Newton's method applied to low-degree polynomials in one complex variable as a global dynamical system on the Riemann sphere. His idea opened an immense field of mathematical research that we now call holomorphic dynamics. Roughly speaking, the idea behind the dynamical system approach to find the solutions of a given equation is to {\it construct} a function (also called method)  with two basic properties. On the one hand, the solutions  of the non-linear equation are (super)-attracting fixed points of the method, so initial conditions close to the solutions converge to the solutions as fast as possible. On the other hand, it is sought that for most initial conditions the iterations generated by the function converge to those solutions of the equation.

These dynamical systems lie at the intersection of theoretical work on dynamical systems and applied mathematics. The rationale supporting this claim is that, by studying the topological and geometrical properties of the (immediate) basins of attraction associated with the fixed points of the dynamical system—a natural source of theoretical inquiry—we deepen our understanding of how to design better methods and how to identify good sets of initial conditions to explicitly find (or approximate) the desired solutions of non-linear equations.

Nonetheless, the general problem of solving non-linear equations requires some concreteness: regularity of the equation, finiteness or infiniteness of solutions, dimension and topology of the phase space, etc. For each case, the methodology and output(s) might differ substantially. Accordingly,  we introduce next some notation and describe the framework in which we will develop our study.  This introduction is expanded in Section \ref{section:prelim}, where we present the main concepts and results of complex dynamics used along the paper. For a more detailed introduction to the theory of rational iteration see, for instance, \cite{CarlesonGamelinBook, BeardonBook, MilnorBook}.

 Let $p$ be a  polynomial in one complex variable. Then, by applying a root-finding algorithm to $p$ we obtain a (rational) map $F_p:\hat{\mathbb C}\to \hat{\mathbb C}$, where $\hat{\mathbb C}:=\mathbb C \cup \{\infty\}$ denotes the Riemann sphere. Moreover, every root $\alpha$ of $p$ is a (locally) attracting fixed point of $F_p$: $F_p(\alpha)=\alpha$ and if an initial condition $z_0$ is close enough to $\alpha$ then the sequence of iterates $\{F_p^n(z_0)\}$ converges to $\alpha$.

 When we study the dynamical system induced by the iterates of  $F_p$ we aim to understand the asymptotic behaviour of the sequences 
$\{z_n=F_p^n(z_0)\}$ for any seed $z_0\in \mathbb C$. As we already mentioned, if $p(\alpha)=0$ then $F_p(\alpha)=\alpha$  and, since $\alpha$ is locally attracting, we have $\displaystyle \lim_{n \to \infty} z_n= \alpha$ if $z_0\approx \alpha$.  This fact allows to introduce the notion of basin of attraction for any root $\alpha$ of $p$ as
$$
A_{F_p}(\alpha)=\{z\in \mathbb C \ | \ F_p^n(z) \to \alpha \text{ as } n \to \infty \}.
$$
 The immediate basin of attraction of $\alpha$, denoted by $A^{\star}_{F_p}(\alpha)$, is the connected component of $A_{F_p}(\alpha)$ that contains $\alpha$. Notice that  $A_{F_p}(\alpha)$ and   $A^{\star}_{F_p}(\alpha)$ are open sets. The connected components of $A_{F_p}(\alpha)$ are part of the well-known Fatou set of $F_p$  while their boundary, denoted by $\partial A_{F_p}(\alpha)$, belongs to the Julia set (see details in Section \ref{section:prelim}).

The most well-known and studied root-finding algorithm is induced by the map
\begin{equation} \label{eq:Newton}
F_p(z):=N_p(z)=z-\frac{p(z)}{p^{\prime}(z)}
\end{equation} 
known as the {\it Newton's map}. It is well known that  $N_p(\alpha)=\alpha$, $\ \alpha \in \mathbb C,$ if and only if $p(\alpha)=0$. Using appropriate coordinates, we also have $N_p(\infty)=\infty$.
  Moreover, if $\alpha$ is a simple root of $p$, then it is a super-attracting fixed point of $N_p$, i.e.\ $N^{\prime}_p(\alpha)=0$ (multiple roots of $p$ are attracting but not super-attracting, $0<|N^{\prime}_p(\alpha)|<1$). On the other hand, $z=\infty$ is a repelling fixed point,  i.e.\ $|N^{\prime}_p(\infty)|>1$. 

In \cite{HowToNewton} the authors show how to use a deep understanding on the topology and geometry of $A_{N_p}(\alpha)$ in order to generate a universal set $\mathcal S_d$ of initial conditions (only depending on the degree $d$ of the polynomial $p$) such that iterating those points one can find all roots of {\it any degree $d$} polynomial. Their result is based on the fact that $A^{\star}_{N_p}(\alpha)$ are unbounded and simply connected and a key part of the argument is based on the fact that $N_p$ has no other finite fixed points than the roots of $p$. See also \cite{ConnectivityJulia}.

Despite Newton's method plays a central role among all root-finding algorithms it has some limitations. For instance its local order of convergence is two. Over the years there has appeared alternatives root-finding algorithms with higher convergence rate so that, locally, they approach faster the roots of $p$. Among the cubic root-finding algorithms (local order of convergence three) the most well-known are Halley's and Traub's methods: $H_p$ and $T_p$, respectively.  Halley's method belongs to the (discrete) family of {\it König's} root-finding methods (see \cite{BufHen}) with local degree of convergence $n$ to simple roots of a polynomial. For $n=2$ König's method coincides with Newton's method while Halley's method corresponds to $n=3$.  Traub's method is a modification of Newton's method where $p'(z)$ is reused every other iterate (see \cite{Traub_Book}). Halley's method $H_p$ and Traub's method $T_p$ are obtained as follows:
$$
H_{p}(z):=z-\frac{p(z)p^{\prime}(z)}{\left(p^{\prime}(z)\right)^2-\frac{1}{2}p(z)p^{\prime\prime}(z)} \quad\quad  \mbox{and} \quad\quad T_p (z) :=  N_p(z)  - \frac{ p(N_p(z))}{p'(z)}.  
$$ 

As mentioned earlier, both maps define cubic root-finding algorithms that offer better local convergence to simple roots than Newton’s method. However, this local improvement comes at a cost when we consider the associated global dynamics on $\hat{\mathbb{C}}$. We previously argued that a key property of Newton’s method is that its finite fixed points coincide with the roots of the polynomial $p$.  This property, which is key for the understanding of the dynamics of $N_p$, no longer holds neither  for Halley's nor for Traub's method.

The other key property used in  \cite{HowToNewton} to obtain the universal set of initial conditions for Newton's method $N_p$ is the fact that the immediate basins of attraction of roots are simply connected and unbounded under $N_p$. In this paper we propose partial answers to the natural question of whether the immediate basins of attraction of the roots of $p$ are unbounded and simply connected under Halley's and Traub's methods. We already know partial results about this matter. In \cite{CEGJ2023} it is shown that the dynamical plane might have bounded immediate basins of attraction for Traub's method applied to cubic polynomials but those fixed points do not correspond to zeros of the polynomial. In the same paper the authors show the unboundedness and simple connectivity of the immediate basins of attractions of the roots of $p$ for Traub's method applied to a concrete family of (symmetric)  degree $d$ polynomials.  Moreover, the authors also conjectured that these properties hold for $T_p$ for every polynomial $p$. The study of the unboundedness and simple connectivity of basins of attraction of roots under root-finding algorithms is indeed a very active topic of research (see, for instance, \cite{Feliks89,ConnectivityJulia,BFJK18,Par,NP25}.

The main results of this paper are the following. Let $d\geq 2$ and consider family of  degree $d+1$ polynomials $p_d(z)=z(z^d-1)$. Denote by $\alpha_k, \ k=0,\ldots, d-1$ the $d$th-roots of unity and by $A^{\star}_{H_d}(\alpha_k)$ and $A^{\star}_{T_d}(\alpha_k)$ their immediate basins of attraction for Halley's and Traub's methods. Also denote by $A^{\star}_{H_d}(0)$ and $A^{\star}_{T_d}(0)$ the immediate basins of attraction of $z=0$  for Halley's and Traub's methods.

\begin{thmA}
Let $d\geq 2$ and consider $H_d(z)$ to be Halley's map applied to $p_d$. Then, the following statements hold:  
\begin{itemize}
\item[(a)] The sets $A^{\star}_{H_d}(\alpha_k), \ k=0,\ldots,d-1,$ as well as $A^{\star}_{H_d}(0)$  are simply connected.
\item[(b)] The sets $A^{\star}_{H_d}(\alpha_k), k=0,\ldots,d-1,$  are unbounded.
\item[(c)] For $d=2,3,4$ the set $A^{\star}_{H_d}(0)$ is unbounded.
\item[(d)] For $d \geq 5$ the set $A^{\star}_{H_d}(0)$ is bounded.
\end{itemize}
\end{thmA}

\begin{thmB}
Let $d\geq 2$ and consider $T_d(z)$ to be Traub's map applied to $p_d$. Then,  $A^{\star}_{T_d}(\alpha_k),\ k=0,\ldots, d-1$ as well as $A^{\star}_{T_d}(0)$ are simply connected and unbounded.
\end{thmB}

In light of Theorem B, the main conjecture in \cite{CEGJ2023} remains open. In contrast, Theorem A shows that the unboundedness of immediate basins is not a general property of Halley’s method. In fact, our counterexample was inspired by Proposition 6 in \cite{BufHen}.

The paper is organized as follows. In Section \ref{section:prelim} we provide the necessary background to support both the technical and main results presented. In Section \ref{section:TheoremA} we prove Theorem A and present additional findings on Halley’s method applied to the family $p_d$. Finally, in Section \ref{section:TheoremB} we prove Theorem B.

{\it Acknowledgments}. We thank  Christian Henriksen for very helpful discussions.  

\section{Preliminaries}\label{section:prelim}

We start this section introducing basic concepts of complex rational iteration. See \cite{CarlesonGamelinBook, BeardonBook, MilnorBook} for a more detailed overview of the field of complex dynamics and proofs of the results mentioned below. Let $R:\hat{\mathbb C} \to \hat{\mathbb C}$ be a rational (analytic) map of degree $d$. We say that $z\in \hat{\mathbb C}$ (using appropriate coordinates if $z=\infty$) is a point in the {\it Fatou set} of $R$, $\mathcal F(R)$, if there exist a sufficiently small neighbourhood $U$ of $z$ such that the family of iterates $\{R^n|_U\}$ is normal. Otherwise we say that $z$ belongs to the {\it Julia set} of $R$, $\mathcal J(R)$.
 By definition, these sets provide a partition of $\hat{\mathbb C}$ and it is known they are dynamically invariant. Moreover, $\mathcal J(R)$ is closed and $\mathcal F(R)$ is open. Each connected component of $\mathcal F(R)$ is called a {\it Fatou component} or {\it Fatou domain}. Since $\mathcal F(R)$ is totally invariant, Fatou components are mapped among themselves under iteration of $R$. By Sullivan's No-Wandering Domain Theorem \cite{Su}, every Fatou component is periodic or eventually periodic. For instance, if $z_0$ is an attracting fixed point for $R$, that is, $R(z_0)=z_0$ and $|R^{\prime}(z_0)|<1$, then the immediate attracting basin of $z_0$, $A_{R}(z_0)$, is a periodic (fixed) Fatou component and its iterated preimages are eventually periodic Fatou components. It follows from the Classification Theorem that every periodic Fatou component of a rational map belongs to the immediate basin of attraction of an attracting or parabolic fixed point, to a simply connected rotation domain (Siegel disk), or to a doubly connected rotation domain (Herman ring).
 Finally,   we say that $z_0$ is a {\it critical point} if $R^{\prime}(z_0)=0$. The  image $R(z_0)$ of a critical point is called a {\it critical value}. It is well known that critical points and critical values play a central role in holomorphic dynamics. For instance, every periodic Fatou component can be related to a critical point: attracting and parabolic basins contain, at least, a critical point while the orbit of a critical point accumulates on the boundary of every rotation domain. 

 Next we introduce some results that are used along the paper. The following result is proven in \cite{BFJK18}.  It uses the concept of weakly repelling fixed point: we say that $z_0$ is a {\it weakly repelling fixed point} for $R$ if $R(z_0)=z_0$ and $|R^{\prime}(z_0)|>1$ or $R^{\prime}(z_0)=1$. 

\begin{lemma}[\bf Boundary maps out] \label{lem:mapout}
Let $\Omega \subset \mathbb C$ be a bounded domain with finite Euler characteristic and let $f$ be a meromorphic map 
in a neighbourhood of $\overline{\Omega}$. Assume that there exists a component
$D$ of \ $\hat{\mathbb C} \setminus f(\partial\Omega)$, such that:
\begin{itemize}
\item[$(a)$] $\overline{\Omega} \subset D$,
\item[$(b)$] there exists $z_0 \in \Omega$ such that $f(z_0) \in D$. 
\end{itemize}
Then, $f$ has a weakly repelling fixed point in $\Omega$. 
\end{lemma}

Let $U,V\in \mathbb C$ open domains. We say that $f:U\to V$ is a degree $d$ proper map if every point in $V$ has $d$ preimages (counting multiplicity) and the preimage of every compact set $K\in V$ is compact in $U$. For instance if $V\subset \mathbb C$ and $U$ is a connected component of $R^{-1}(V)$ then $R:U\to V$ is a proper map. We recall the Riemann-Hurwitz formula (see for instance \cite{Ste}), which we use in order to study connectivities of Fatou components.

\begin{theorem}{(Riemann-Hurwitz formula)}
	\label{theorm:RH}
	Let $U,V \subset \hat{\mathbb{C}}$ be two connected domains of connectivity $m_U, m_V \in \mathbb{N}^{*}$ and let $f: U \to V$ be a degree $d$ proper map branched over $r$ critical points, counted with multiplicity. Then, 
	$$m_U-2=d(m_V-2)+r .$$
\end{theorem}

We also use the following corollary of the Riemann-Hurwitz formula (compare \cite[Corollary 2.2]{CFG1}).
\begin{corollary}
	\label{corollary:HR_coro}
	Let $U\subset \hat{\mathbb C}$ be an open set and let $f: U \to f(U)$ be a proper holomorphic map. Then, the following statements hold:
	\begin{enumerate}
		\item[(a)] If $f(U)$ is doubly connected and $f$ has no critical points in $U$, then $U$ is doubly connected.
		\item[(b)] If $f(U)$ is simply connected and $f$ has at most one critical point in $U$ (not counting multiplicities), then $U$ is simply connected. 
	\end{enumerate}
\end{corollary}

\section{Proof of Theorem A} \label{section:TheoremA}

Let $p$ be a polynomial of degree $d$ with all roots being simple. Halley's map (see  \cite{BufHen}) applied to $p$ is given by the formula
\begin{equation}\label{eq:halley}
H_{p}(z):=z-\frac{p(z)p^{\prime}(z)}{\left(p^{\prime}(z)\right)^2-\frac{1}{2}p(z)p^{\prime\prime}(z)}
\end{equation}
One can easily check that the roots of $p$ correspond to superattracting fixed points of  {\bf $H_p$} while the zeros of $p^{\prime}$ are repelling fixed points of $H_p$ with multiplier 3.

In this section we prove Theorem A, that is,  we study the boundedness and simple connectivity of the immediate basins of attraction  under the map $H_d$ obtained applying  Halley's method to the degree $d\geq 2$  family of polynomials given by 
\begin{equation} \label{eq:Ham_pd}
p_d(z):=z(z^d-1). 
\end{equation}
We denote by $\{0,\alpha_0=1,\alpha_1,\ldots,\alpha_{d-1}\}$ the $d+1$ roots of $p_d$, where
$$
\alpha_k=\exp\left(\frac{2k\pi i}{d-1}\right), \ k=0,\ldots d-1.
$$

To simplify notation we denote by $A_d(0):=A_{H_d}(0)$ and $A_d(\alpha_k):=A_{H_d}(\alpha_k)$, $k=0, \ldots, d-1$ the basin of attraction of each of the roots $p$ and  by $A_d^\star\left(0\right):=A^{\star}_{H_d}(\alpha_k)$ and $A_d^\star\left(\alpha_k\right):=A^{\star}_{H_d}(\alpha_k),\ k=0,\ldots d-1$, the corresponding immediate basins of attraction. We start by proving some (rather immediate) lemmas we will use later to prove Theorem A. Some of them are proved (or they are direct consequences)  in \cite{BufHen}, but we state and prove them here for completeness.

\begin{lemma}[\bf Dynamical symmetry] \label{lem:technical1}
 Consider the family of polynomials given in \eqref{eq:Ham_pd} and define the two auxiliary functions
 \begin{equation} \label{eq:q_Q}
 \begin{split}
 &q_d(w)=\frac{dw[(d-1)+(d+1)w]}{2+(d+1)(d-4)w+(d+1)(d+2)w^2}  \\
&Q_d(w)=2(d-1)+(d+1)(d+2)w. 
 \end{split}
 \end{equation}
 The following statements about $\Hfam$ hold.
 \begin{enumerate} [label=(\mbox{\alph*})]
 \item The map $\Hfam$ is given by
 \begin{equation}\label{eq:halley_pd}
 \Hfam(z)=\frac{dz^{d+1}[(d-1)+(d+1)z^{d}]}{2+(d+1)(d-4)z^{d}+(d+1)(d+2)z^{2d}}=:zq_d(z^{d}) \  .
 \end{equation}
In particular, $\Hfam$ has degree $2d+1$. 
\item The derivative of $\Hfam$ is given by 
\begin{equation} \label{eq:deriv_halley_pd}
\Hfam^{\prime}(z)=q_d(z^{d})+dz^{d}q_d^\prime(z^{d})=\frac{d(d+1)z^d(z^d-1)^2}{(2+(d+1)(d-4)z^{d}+(d+1)(d+2)z^{2d})^2}Q_d(z^d)
\end{equation}
\item The point $z=0$ is a  super-attracting fixed point with local degree $d+1$ and, therefore, a critical point of multiplicity $d$. 
\item The $d$th-roots of unity are  super-attracting fixed points with local degree 3 and, therefore, are critical points of multiplicity, at least, 2.
\item The dynamical plane for $\Hfam$ is symmetric with respect to rotation by a $d$th-roots of unity. That is, if $\eta\in \mathbb C$ is such that $\eta^{d}-1=0$ and $\phi(z)=\eta z$, then
$$
\Hfam(z)=\left(\phi\circ \Hfam \circ \phi^{-1}\right)(z).
$$
\end{enumerate}
\end{lemma}

\begin{proof}
Statement (a) is direct by substituting $p$ by $p_d$ in \eqref{eq:halley}. Statements (b) follows from statement (a) and simplifications. Statements (c) and (d) are direct. Statement (e) is equivalent to show that if $\eta^{d}-1=0$ then $\eta\Hfam(z)=\Hfam(\eta z)$,  which is a direct computation from \eqref{eq:halley_pd}. 
\end{proof}

\begin{figure}
	
\centerline{\includegraphics[width=0.7\textwidth]{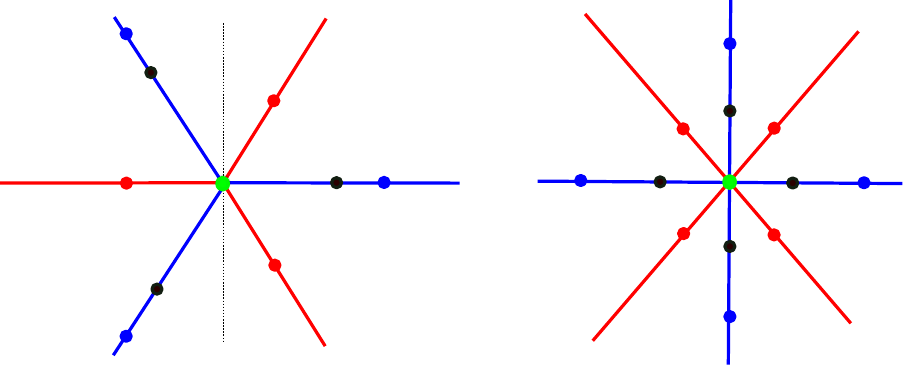}
 \put(-213,85){$c_0$}
 \put(-272,53){ $c_1$}
  \put(-213,35){$c_2$}
   \put(-182,53){ $\alpha_0$}
   \put(-200,52){ $\zeta_0$}
    \put(-268,90){ $\zeta_1$}
   \put(-283,110){ $\alpha_1$}
    \put(-268,24){$\zeta_2$}
   \put(-281,10){$\alpha_2$}
   \put(-13,54){$\alpha_0$}
      \put(-42,53){ $\zeta_0$}
 \put(-91,53){$\zeta_2$}
      \put(-115,55){$\alpha_2$}
    \put(-50,73){ $c_0$}
        \put(-45,42){ $c_3$}
    \put(-80,73){$c_1$}
        \put(-91,43){ $c_2$}
   \put(-60,87){ $\zeta_1$}
 \put(-57,110){$\alpha_1$}
 \put(-61,11){ $\alpha_3$}
 \put(-62,33){ $\zeta_3$}
 }
\caption{\small Distribution of the free critical points and the invariant straight lines. Blue points represent the roots $\alpha_k$  different from zero, red points represent the critical points $c_k$ different from zero, and black points represent the fixed points $\zeta_k$. We sketch the case $d=3$ (left) and the case $d=4$ (right).}
\label{fig:rays1}
\end{figure}
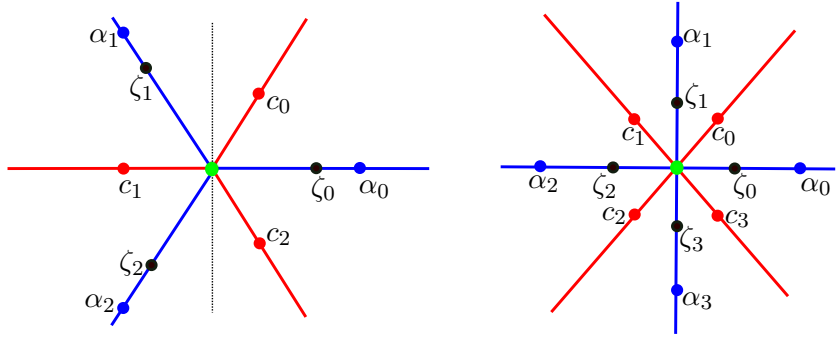

The above lemma explains how the intrinsic symmetry of the family $p_d$ induces a symmetry on the dynamical plane of $H_{p_d}$. In the next lemma, we explore how the symmetry applies to the location and dynamical behaviour of the critical points. In Figure \ref{fig:rays1} we summarize  this structure for $d=3$ and $d=4$.

\begin{lemma}[\bf Critical points and invariant lines] \label{lem:technical2}
The following statements hold.
\begin{enumerate}[label=(\mbox{\alph*})]
\item The critical points of $H_{p_d}$ are $z=0$ (with multiplicity $d$), $z=\alpha_k$ with $ k=0,\ldots d-1$ (with multiplicity 2) and $d$ simple critical points given by 
\begin{equation} \label{eq:critics}
c_{k,d}:=c_k=\left(\frac{2(d-1)}{(d+1)(d+2)}\right)^{\frac{1}{d}} \exp\left(\frac{(2k+1)\pi i}{d}\right), \ k=0,\dots, d-1.
\end{equation}
\item Let $t>0$ and consider the semi-lines 
$$
r_{\ell,d}:=r_\ell=t\exp\left(\frac{\ell\pi i}{d}\right),\ \ell=0,\ldots, 2d-1.
$$
\begin{itemize}
\item[(b.1)] If $\ell$ is even the semi-line $r_\ell$ is forward invariant and  it contains the root $\alpha_{\ell/2}$.

\item[(b.2)] If $\ell$ is odd the semi-line $r_{\ell}$ verifies
\begin{equation}\label{eq:ell_ell_prime}
\Hfam(r_\ell) \subset r_\ell \cup r_{\ell^\prime} \cup \{0\}:=\tilde{r}_\ell, \ \ \mbox{\rm with} \ \  \ell^{\prime} = (\ell+d)\ \mbox{\rm mod}(2d).
\end{equation}
Moreover it contains the critical point $c_{(\ell-1)/2}$.
\end{itemize}
\item Let $t\in\mathbb{R}$ and consider the lines
$$
\tilde{r}_\ell=t\exp\left(\frac{\ell\pi i}{d}\right)=:t\exp \left(i\theta_{\ell}\right),\ \ell=0,\ldots, d-1.
$$
Then, the anticonformal map $\iota_{\ell}(z)=\exp\left(2i\theta_{\ell}\right)\overline{z}$ is a reflection with respect to the line $\tilde{r}_\ell$ that conjugates $\Hfam$ with itself (i.e.\ $\Hfam=\iota_{\ell}\circ \Hfam \circ\iota_{\ell}$).
\end{enumerate}
\end{lemma}

\begin{proof}[Proof of Lemma \ref{lem:technical2}]
The expression of $\Hfam^{\prime}$ in \eqref{eq:deriv_halley_pd}  and \eqref{eq:q_Q} implies statement (a).  To prove statement (b) we just notice that from \eqref{eq:q_Q} and \eqref{eq:halley_pd} we have

\begin{equation} \label{eq:Hfam_restricted}
\Hfam\left(t\exp\left(\frac{\ell\pi i}{d}\right)\right)=\exp\left(\frac{\ell\pi i}{d}\right)\ t\ q_d\left(t^d (-1)^\ell\right). 
\end{equation}

For $\ell$ even we claim that that $\Hfam(r_\ell)= r_\ell$. Indeed, we have $\Hfam(0)=0$, $q_d\left(t^d \right)>0$ and
$$
\lim\limits_{t\to \infty} t\ q_d\left(t^d \right) = \lim\limits_{t\to \infty} \frac{d(d+1)}{(d+1)(d+2)}t =\infty. 
$$

For $\ell$ odd the numerator and denominator of $q_d(-t^d)$ might be positive or negative depending on the value of $d$. Hence \eqref{eq:ell_ell_prime} is satisfied. In fact, some computations show that  for $d<7$ we have that 
$\Hfam(r_\ell) \subsetneqq \tilde{r}_\ell$ while for $d\geq 7$ the denominator of $q_d$ vanishes and then $\Hfam(r_\ell) = \tilde{r}_\ell$ 
  
Statement (c) follows directly from the Schwarz Reflection Principle  using that the lines $\tilde{r}_\ell$ are forward invariant and their points are preserved by the anticonformal map $\iota_{\ell}$. In particular, if $z\in \tilde{r}_\ell$ then $\iota_\ell(z)=z$.
\end{proof}

\begin{lemma}[\bf Fixed points of $\Hfam$] \label{lem:technical3}
The map $\Hfam$ has exactly $2d+2$ fixed points in $\hat{\mathbb C}$. Moreover, the following statements hold.

\begin{enumerate}[label=(\mbox{\alph*})]
\item The point $z=\infty$ is a repelling fixed point with  positive real multipler. 
\item There are $d+1$ superattracting fixed points corresponding to the $d+1$ roots of $p_d$.
\item There are $d$ repelling fixed points corresponding to the $d$ roots of $p_d^{\prime}(z)=(d+1)z^{d}-1$. More precisely, for every $k \in \{0,\ldots ,d-1\}$ there exists a unique  $\zeta_{k,d}\in r_{2k}$  such that  $\Hfam(\zeta_{k,d})=\zeta_{k,d}$  and
$$
\zeta_{k,d} \in \partial A_d^{\star}(0) \cap \partial A_d^{\star}(\alpha_k).
$$
\item Each $\zeta_{k,d}$ is accessible from $A_d^{\star}(0)$.
\end{enumerate}
\end{lemma}

\begin{proof}
Since $\Hfam$ has degree $2d+1$, it has $2d+2$ fixed points. From \eqref{eq:halley_pd} it is clear that $z=\infty$ is a fixed point and standard computations show that the multiplier is
$1+2/d>1$.  This proves statement (a). From \eqref{eq:halley}, the finite fixed points of $\Hfam$ are either zeros of $p_d$, that is $\{0,\alpha_0,\ldots,\alpha_{d-1}\}$, or  zeros of  $p^{\prime}_d$. So, statement (b) follows. 

To prove statement (c) observe that the zeros of $p_d^{\prime}$ denoted by $\zeta_{k}:=\zeta_{k,d}$  are given by 
\begin{equation}\label{eq:fixed_points_pprime}
\zeta_k=\left(\frac{1}{d+1}\right)^{\frac{1}{d}}\exp\left(\frac{2k\pi i}{d}\right) \in r_{2k}, \ k=0,\ldots d-1.
\end{equation}     
Notice that $|\zeta_k|<1$ for all $k$. It remains to show that 
$$
\zeta_k \in \partial A_d^{\star}(0) \cap \partial A_d^{\star}(\alpha_k).
$$

From Lemma \ref{lem:technical1} (dynamical symmetries) it is enough to see that the result is true for $\zeta_0\in (0,1)$. The map $\Hfam$ is a real map and  from, Lemma~\ref {lem:technical1}, the points  $z=0$ and $z=1$ are super-attracting fixed points with local degree $d+1$ and 3, respectively. Moreover, by \eqref{eq:deriv_halley_pd} and \eqref{eq:fixed_points_pprime}  we conclude $\Hfam$  is strictly increasing in the interval $(0, 1)$ with a unique (repelling) fixed point at $\zeta_0 \in (0,1)$. Moreover, we have that $[0,\zeta_0) \subset A_d^\star(0)$ and $(\zeta_0,1] \subset A_d^\star(1)$. Therefore,  $[0,\zeta_0)$ is an invariant access to $\zeta_0$.

From the previous arguments it is clear that all $\zeta_k$'s are accessible from $A_d^{\star}(0)$ and $[0,\zeta_k)\subset r_{2k}$ is an invariant access to $\zeta_k$, $\ k\in\{0,1,\ldots,d-1\}$ (see Definition 2.2 in \cite{bfjkaccesses} for the definition of {\it invariant access}).
\end{proof}

As we will see, the  boundedness of the immediate basin of attraction of $z=0$, that is $A_d^{\star}(0)$, depends on the degree $d$ of the polynomial $p_d$. Nonetheless, this is not the case for the rest of immediate basins of attractions.  

\begin{prop} \label{prop:halley_alphak_unbounded}
 Let $d\geq 2$. Then, for each $k=0,\ldots, d-1$ the set $A_d^{\star}(\alpha_k)$ is unbounded.
\end{prop}

\begin{proof}
To prove the proposition we notice that from Lemma \ref{lem:technical1} it is enough to show that $A_d^{\star}(1)$ is unbounded.  Since  the real line is invariant under $\Hfam$,  $z=1$ is a superattracting fixed point of  $\Hfam$, there are no fixed points in $(1,\infty)$, and  $\Hfam$ restricted to $(1,\infty)$ is strictly increasing (its derivative is positive from \eqref{eq:deriv_halley_pd}), we have that interval $(1,\infty) \in A_d^\star(1)$. By symmetry, we conclude that the basins of attraction $A_d^\star(\alpha_k)$ are unbounded for each $ k=0,\ldots,  d-1$.
\end{proof} 

Before discussing the boundedness of $A_d^\star(0)$
depending on the degree of the polynomial $p_d$, we prove the simple connectivity of the immediate basins of attraction of all the roots of $p_d$.

\begin{prop}\label{prop:simple_connected_Halley}
 The  immediate basins $A_d^\star\left(0\right)$ and $A_d^\star\left(\alpha_k\right), k=0,\ldots, d-1$, are simply connected. 
 \end{prop}
 
 \begin{proof}
We argue by contradiction. Let $A_d^\star(\beta)$ be a multiple connected Fatou component,  where $\beta \in \{0,\alpha_0,\ldots \alpha_{d-1}\}$. Consider a simply connected set $U \subset A_d^\star(\beta) $ such that $\beta \in U$ and $\overline{\Hfam(U)} \subset U$ (this is in fact an absorbing domain for $A_d^\star(\beta)$ and its existence is a consequence of the local B\"otcher coordinates around $z=\beta$). Set $U_\ell:=\Hfam^{-\ell}(U), \ \ell>0$  be the connected component of $\Hfam^{-\ell}(U)$ containing 
 $\beta$. If $A_d^\star(\beta)$ is multiply connected there must exists $\ell_0>0$ such that $U_{\ell},\ \ell=0,\ldots, \ell_0$, are simply connected but $U_{\ell_0+1}$ is multiply connected. Notice that $\partial U_\ell, \ \ell>0,$ belongs to $A_d^\star(\beta)$. By construction, $\hat{\mathbb C}\setminus  U_{\ell_0+1}$ must have at least one bounded connected component, say $\Omega$, such that $\Hfam(\partial \Omega)=\partial U_{\ell_0}$. Moreover, $\Hfam(\Omega)$ must cover $\hat{\mathbb C}\setminus  U_{\ell_0}$ since, otherwise, $\Hfam(\Omega)\subset U_{\ell_0}$ and so $\Omega \subset A_d^\star(\beta)$, which is a contradiction.  In particular, $\Omega$ contains a pole. 
 
Denote by $D$ the unbounded connected component of  
$\hat{\mathbb C}\setminus  \partial \Hfam(\Omega)$, that is, $D:=\hat{\mathbb C}\setminus \overline{U_{\ell_0}}$. The sets $D$, $\Omega$, and $\partial \Omega$ are under the hypothesis of Lemma \ref{lem:mapout} since there is a pole $p \in \Omega$ (playing the role of $z_0$ in Lemma  \ref{lem:mapout}). Consequently, by Lemma  \ref{lem:mapout} we have that $\Omega$ contains a weakly repelling fixed point, say $\zeta \in \Omega$. From Lemma \ref{lem:technical3} we have $$
\zeta=\zeta_{k_0} \in \partial A_d^{\star}(\alpha_{k_0})
$$ 
for some ${k_0}\in\{0,\ldots, d-1\}$. From Proposition \ref{prop:halley_alphak_unbounded} we have that the sets $A_d^{\star}(\alpha_k),\ k=0,\ldots, d-1$, are  all unbounded.  Using the unboundedness of   $A_d^{\star}(\alpha_{k_0})$,  if $\beta\neq \alpha_{k_0}$ we have $\partial \Omega \cap   A_d^{\star}(\alpha_{k_0}) \ne \emptyset$, a contraction since $\partial \Omega \subset A_d^\star(\beta)$. If  $\beta=\alpha_{k_0}$, since $\zeta_{k_0}\in \partial A_d^\star(0)$, we conclude that $A_d^\star(0)\subset \Omega$. In particular, $\zeta_k\in\Omega$ for all $k\in\{0,\ldots, d-1\}$, which leads to a contradiction as before.

 \end{proof}
 
 The following proposition is needed to discuss the boundedness of $A_d^\star(0)$ in terms of $d$.
 
\begin{prop}\label{prop:number_critics}
The following statements hold.
\begin{itemize}
\item[(a)]  If $d\geq 5$ then $\Hfam: A_d^\star(0) \mapsto A_d^\star(0)$ has degree $d+1$ and the only critical point in $A_d^\star(0)$ is $0$ (with multiplicity $d$).
\item[(b)]  If $d=2,3,4$, then   $\Hfam: A_d^\star(0) \mapsto A_d^\star(0)$ has degree $2d+1$ and it contains the simple critical points  $c_k$ for $k=0, \ldots, d-1$ and the critical point 0 with multiplicity $d$.
\end{itemize}
 \end{prop}

\begin{proof}
In the proof we abuse notation not showing the dependence of the critical in $d$. From the above proposition we know that $A_d^*(0)$ is simply connected. Using the Riemann-Hurwitz formula (Theorem \ref{theorm:RH}) we can  obtain the degree of  $\Hfam: A_d^\star(0) \mapsto A_d^\star(0)$  only computing the number of critical points inside $A_d^\star(0)$. We write $c_{k,d}=c_k$ if there is no confusion.

We first compute $v_{k,d}=v_k=\Hfam(c_k),\ k=0,\ldots, d-1$. From \eqref{eq:critics},  we can write  
$$
 c_k = a_d\exp\left(\frac{(2k+1)\pi i}{d}\right) ,\quad \mbox{with} \quad a_d=\left(\frac{2(d-1)}{(d+1)(d+2)}\right)^{\frac{1}{d}}>0.
$$
Then, from \eqref{eq:halley} and some computations we have
$$
v_k=\Hfam(c_k)=c_k \left(\frac{-da_d^d((d-1)-(d+1)a_d^d)}{2-(d+1)(d-4)a_d^d+(d+1)(d+2)a_d^d}\right) =c_k\left(\frac{(d-1)^2}{(d+2)(d-7)}\right)
$$

We start proving statement (a). If $d=7$ the critical points  $c_k, \ k=0,\ldots, 6 $,  are also poles and hence they cannot be in  $A_d^{\star}(0)$ ($z=\infty$ is a repelling fixed point).  So, $z=0$ is the only critical point in $A_d^{\star}(0)$ with multiplicity $d$ (see Lemma \ref{lem:technical2}).

Next we study the cases  $d=5,6$ and $d >7$.  First we notice that if $d=5,6$ we have 
\begin{equation}\label{eq:argument_ck1}
1<\left|\frac{(d-1)^2}{(d+2)(d-7)}\right|  \quad \mbox{and} \quad \mbox{Arg}(v_k)=\mbox{Arg}(c_k)+\pi,
\end{equation}
while if $d>7$ we have 
\begin{equation}\label{eq:argument_ck2}
1<\left|\frac{(d-1)^2}{(d+2)(d-7)}\right|  \quad \mbox{and} \quad \mbox{Arg}(v_k)=\mbox{Arg}(c_k).
\end{equation}
Hence, in both cases we obtain 

\begin{equation} \label{eq:v_and_c}
|v_{k}|=|\Hfam(c_{k})|=|c_k|\left|\frac{(d-1)^2}{(d+2)(d-7)}\right|>|c_k| \quad \mbox{and} \quad c_k,v_k \in \tilde{r}_{2k+1}.
\end{equation}
From this we deduce that $c_k \not\in A_d^{\star}(0)$.  Indeed, assume otherwise and take $U$ to be the maximal domain of definition of the  B\"otcher coordinate around $z=0$ (see, for instance, \cite[Theorem~9.3]{MilnorBook}). It can be shown that $U$ is invariant with respect to the symmetries in dynamical plane described in  Lemma \ref{lem:technical1} and Lemma~\ref{lem:technical2}~(c)). By the symmetry described in  Lemma \ref{lem:technical1}, we have that $c_k\in \partial U$ for all $k\in\{0,\cdots, d-1\}$. Let $\theta_{2k+1}=(2k+1)\pi/d$. Since $U$ is simply connected and it is preserved by the anticonformal reflections described in Lemma~\ref{lem:technical2}~(c), it follows that $U\cap \tilde{r}_{2k+1}=\{te^{i\theta_{2k+1}}, 0<t<|c_k|\}$ and, since $|c_k|<|v_k|$, we conclude that $v_k\notin U$. This is a contradiction since $v_k=f(c_k)\in U$ by construction. We can then conclude that $A_d^\star(0)$ does not contain any critical point other than $z=0$ and, hence, the degree of $\Hfam$ restricted to $A_d^\star(0)$ is $d$. This proves (a).

Now we deal with statement (b). We need to show that $c_k\in A_d^{\star}(0),\ k=0,\ldots, d-1$, for $d=2,3,4$.  Notice that it is enough to show that $c_0 \in A_d^{\star}(0)$ since the symmetry described in Lemma  \ref{lem:technical1} would imply $c_k \in A_d^{\star}(0)$ for all possible $k$. We prove the result case by case. 

Let $d=2$.  By construction $c_{0,2}\in r_{1,2} \subset \tilde{r}_{1,2}$  (we add the dependence $c$, $r$, and $\tilde{r}$ in $d$ to clarify the argument). Moreover, from \eqref{eq:Hfam_restricted}, the dynamics over the invariant line is governed by the real function 
$$
B_2(t)=tq_2(-t^2)=-\frac{2t^2(1-3t^2)}{2+6t^2+12t^4}.
$$ 
We claim that $B_2(t)$ is an even function such that $B_2^{\prime}(0)=0$, it has two other critical points at $\pm \tau_{0,2} \approx \pm 0.3558$, it is strictly decreasing in the interval $(0,\tau_{0,2})$, it is strictly increasing in the interval $(-\tau_{0,2},0)$, and 
$[-\tau_{0,2},\tau_{0,2}] \subset (-\tau_{0,2},0)$. Altogether implies that 
$$
L_2:=\left\{t\exp\left(\frac{\pi i}{2}\right),\  t\in[0,\tau_{0,2}]\right\} \subset A_2^{\star}(0), 
$$
and, hence, $c_{0,2} \in A_2^{\star}(0)$. The claims follow from direct computations.

Let $d=3$. By construction $c_{0,3}\in r_{1,3} \subset \tilde{r}_{1,3}$. Moreover, from \eqref{eq:Hfam_restricted}, the dynamics over the invariant line is governed by the real function 
$$
B_3(t)=tq_3(-t^3)=-\frac{2t^3(2-3t^3)}{2+4t^2+20t^6}.
$$ 
In this case we claim that $B_3^{\prime}(0)=0$, there are two critical points $\tilde{\tau}_{0,3}\approx -1.0420$ and $\tau_{0,3}\approx 0.5489$, 
it is strictly decreasing in the interval $(\tilde{\tau}_{0,3},\tau_{0,3})$, and 
$B_3^2([0,\tau_{0,3}]) \subset [0,\tau_{0,3})$. Altogether implies that 
$$
L_3:=\left\{t\exp\left(\frac{\pi i}{3}\right),\  t\in[0,\tau_{0,3}]\right\} \subset A_3^{\star}(0), 
$$
and, hence, $c_{0,3} \in A_3^{\star}(0)$. The claims follow from direct computations.

Let $d=4$ (this case is similar to $d=2$). By construction, $c_{0,4}\in r_{1,4} \subset \tilde{r}_{1,4}$. Moreover, from \eqref{eq:Hfam_restricted}, the dynamics over the invariant line is governed by the real function 
$$
B_4(t)=tq_4(-t^4)=-\frac{2t^4(3-5t^4)}{2+30t^8}.
$$ 
We claim that $B_4(t)$ is an even function such that $B_4^{\prime}(0)=0$, it has two other critical points at $\pm \tau_{0,4} \approx \pm 0.6421$, it is strictly decreasing in the interval $(0,\tau_{0,4})$, it is strictly increasing in the interval $(-\tau_{0,4},0)$, and 
$[-\tau_{0,4},\tau_{0,4}] \subset (-\tau_{0,4},0)$. Altogether implies that 
$$
L_4:=\left\{t\exp\left(\frac{\pi i}{4}\right),\  t\in[0,\tau_{0,4}]\right\} \subset A_4^{\star}(0), 
$$
and, hence, $c_{0,4} \in A_4^{\star}(0)$. The claims follow from direct computations.

\end{proof}

\begin{figure}[htb]
	\centering
	\subfigure[\scriptsize{$d=3$ }]{\begin{tikzpicture}
			\begin{axis}[width=0.4\textwidth, axis equal image, scale only axis,  enlargelimits=false, axis on top]
				\addplot graphics[xmin=-3,xmax=3,ymin=-3,ymax=3] {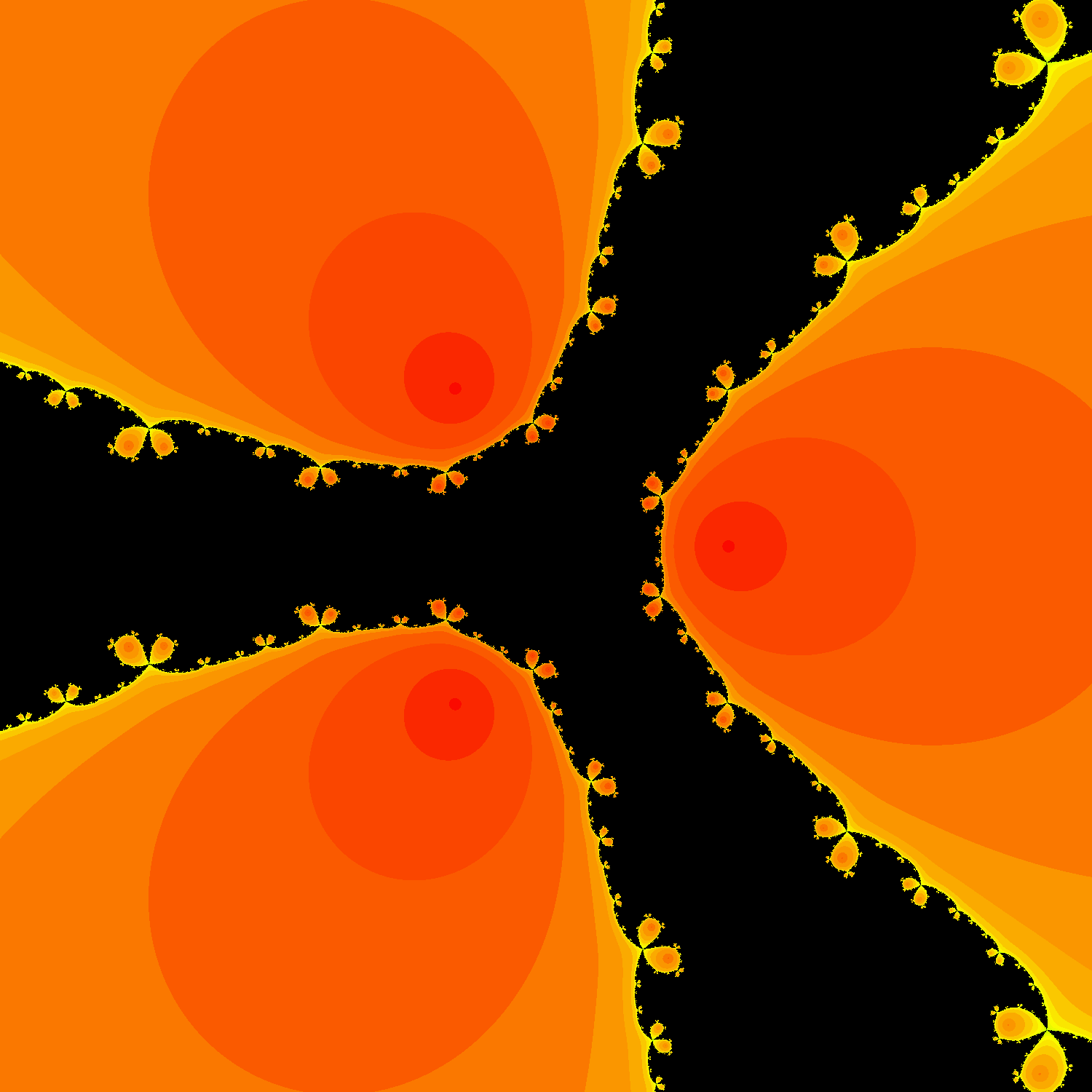};
			\end{axis}
	\end{tikzpicture}}
	\subfigure[\scriptsize{$d=4$ }  ]{	\begin{tikzpicture}
			\begin{axis}[width=0.4\textwidth, axis equal image, scale only axis,  enlargelimits=false, axis on top]
				\addplot graphics[xmin=-3,xmax=3,ymin=-3,ymax=3] {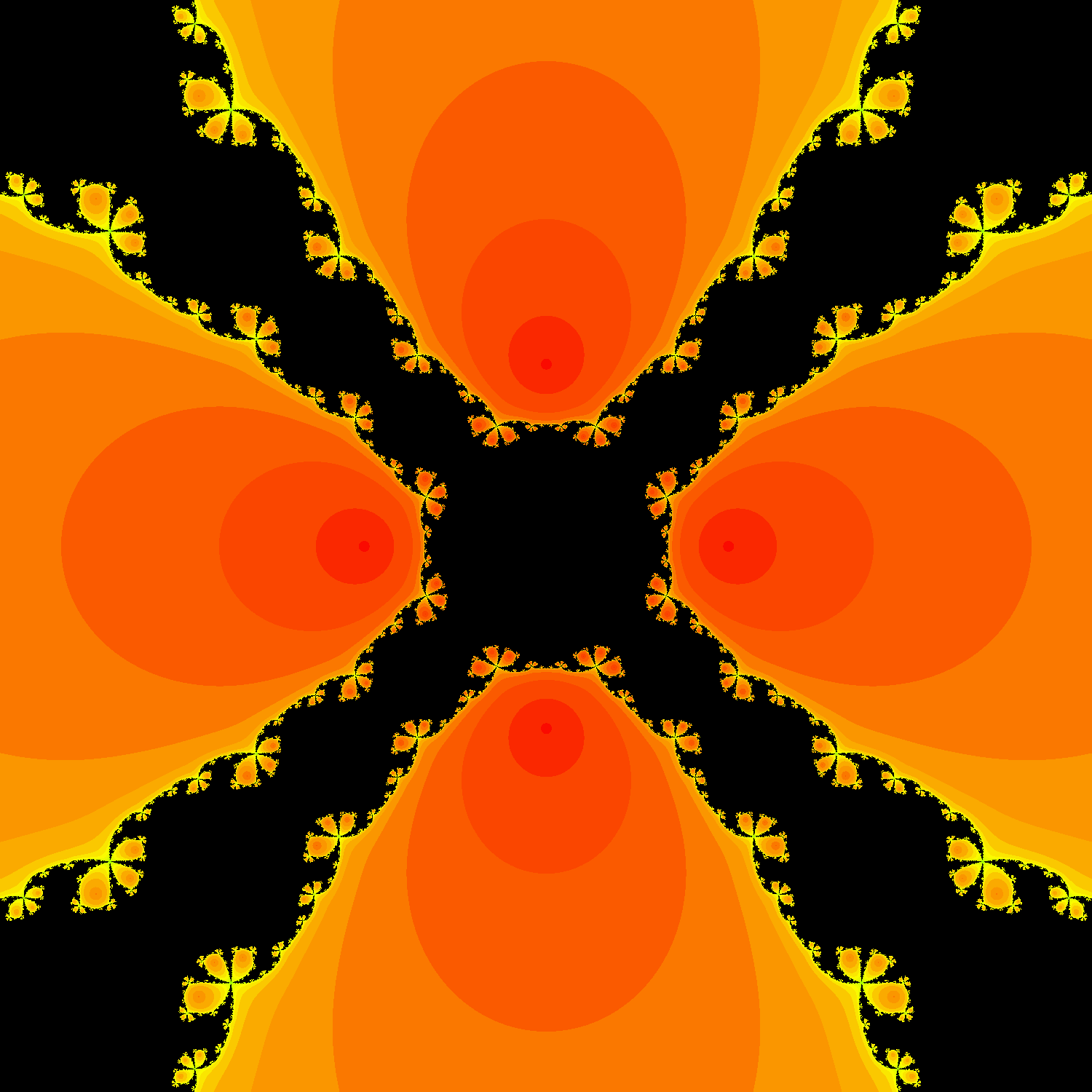};
			\end{axis}
	\end{tikzpicture}}
	\subfigure[\scriptsize{$d=5$}  ]{	\begin{tikzpicture}
			\begin{axis}[width=0.4\textwidth, axis equal image, scale only axis,  enlargelimits=false, axis on top]
				\addplot graphics[xmin=-3,xmax=3,ymin=-3,ymax=3] {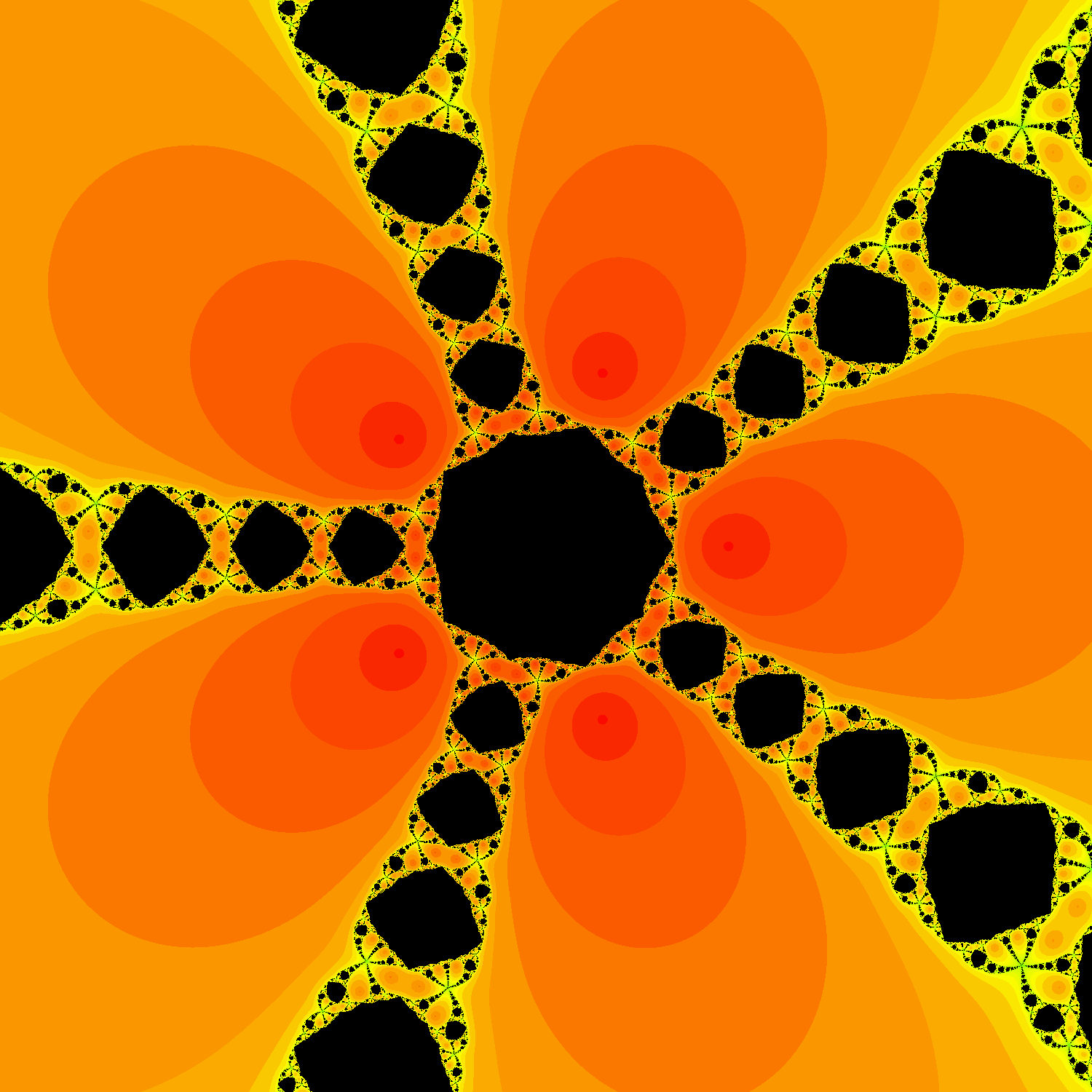};
			\end{axis}
	\end{tikzpicture}}
	\subfigure[\scriptsize{$d=6$}  ]{ \begin{tikzpicture}
			\begin{axis}[width=0.4\textwidth, axis equal image, scale only axis,  enlargelimits=false, axis on top]
				\addplot graphics[xmin=-3,xmax=3,ymin=-3,ymax=3] {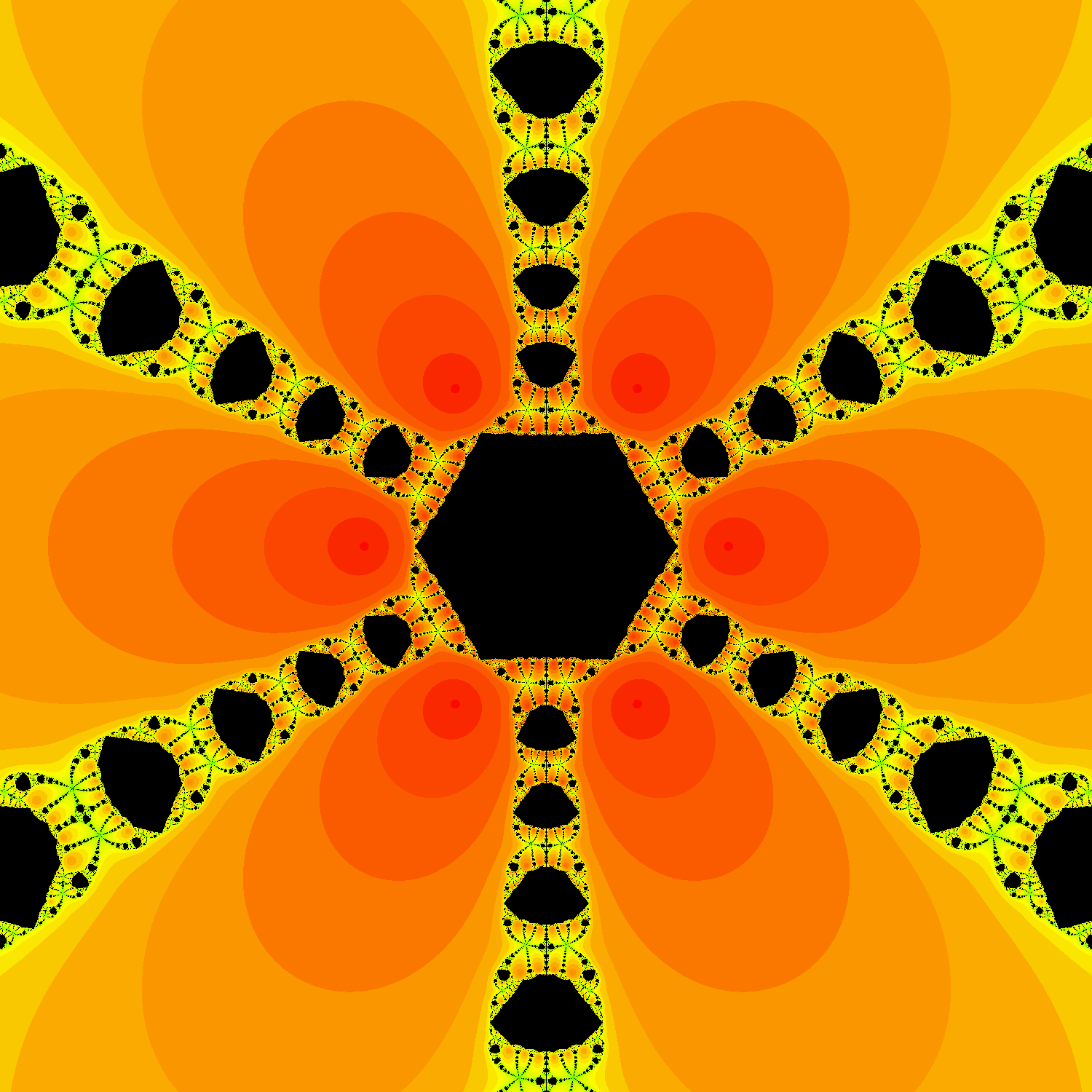};
			\end{axis}
	\end{tikzpicture}}
	\caption{\small{Dynamical planes of Halley's method applied to the family of polynomials $p_d(z)=z(z^d-1)$. We see in black the basin of attraction of 0. The scaling of colours is used to plot the basins of attraction of the $d$th-roots of the unity.}}
	\label{fig:HalleyD}
\end{figure}

\begin{theorem}
$A_d^\star(0)$ is a bounded set if and only if $d\geq 5$. 
\end{theorem}

\begin{proof}
Assume first $d\geq 5$. From Proposition \ref{prop:number_critics}(a) the only critical point of $A_d^\star(0)$ is $z=0$ and  the degree of $\Hfam$ restricted to $A_d^\star(0)$ would be  $d+1$.  Moreover, from Proposition \ref{prop:simple_connected_Halley}, we have that $A_d^\star(0)$  is simply connected. Since there are no extra critical points on $A_d^\star(0)$ the B\"ottcher (local) coordinate extends globally over the whole immediate basin and the dynamics restricted to $\partial A_d^\star(0)$ is conjugated to the map $\theta \mapsto (d+1)\theta,  \theta \in [0,1)$. Hence there are exactly $d$ (internal) fixed rays landing at $d$ fixed points in $\partial A_d^\star(0)$. We already found these landing rays in Lemma \ref{lem:technical3}(c)-(d) and they land to the points $\zeta_k, \ k=1,\ldots d-1$. So, there are no further fixed ray that might land at infinity, which is a repelling fixed point with positive real multiplier (see Lemma \ref{lem:technical3}(a)), say $\lambda>1$.  To finish the argument we notice that under our setting of global B\"ottcher coordinates if $\infty \in \partial A_d^{\star}(0)$ then $z=\infty$ would be the landing point of a periodic ray (see \cite{MilnorBook,DouHub}). Since we proved no fixed rays might land at $z=\infty$ then the ray(s) must be periodic of period at least two, which is not possible because the dynamics in a sufficiently small neighbourhood of $z=\infty$ is conjugated to    
$z\mapsto \lambda z,\ \lambda \in \mathbb R^+\setminus \{0\}$. Therefore, $\infty \not\in \partial A_d^\star(0)$ or, equivalently, the immediate basin $A_d^\star(0)$ is bounded.

Assume now that $d<5$. From Proposition \ref{prop:number_critics}(b), the degree of $\Hfam$ restricted to $A_d^\star(0)$ is  $2d+1$ (notice that this degree is significantly higher than the case above $d\geq 5$). In particular, since the total degree of $\Hfam$ is also $2d+1$ (see Lemma \ref{lem:technical1}(a)) we have that $A_d^\star(0)$ is totally invariant.

From Proposition \ref{prop:simple_connected_Halley} we have that  $A_d^\star(0)$ is simply connected, and so we might consider $\varphi:\mathbb D \to A_d^\star(0)$ the Riemann map. Normalizing we assume that $\varphi$  sends the interval (0,1) into $\mathbb{R}^+\cap A_d^\star(0)$. The map $\varphi$ can be used to define 
$$
g:=\varphi^{-1}\circ \Hfam \circ \varphi: \mathbb D \to \mathbb D$$ 
the inner function associated to $\Hfam|_{A_d^\star(0)}$. In fact  $g$ extends, using Schwartz's reflexion, to a degree $2d+1$ rational map $G:\hat{\mathbb C} \to \hat{\mathbb C}$. By construction $G$ has $2d+2$ fixed points in $\hat{\mathbb C}$. Those are $z=0$, $z=\infty$, and $\beta_k \in \mathbb S^1$  with $k=1,\ldots,2d$. We notice that $z=0$ is  a simple fixed point of $\Hfam$ since it is a simple solution of $\Hfam(z)-z=0$ 

It is known (see, for instance, \cite{bfjkaccesses}) that from the conjugacy between $g$ and $\Hfam$ through $\varphi$, each of the $\beta_k$'s characterize a different access to a fixed points in $\partial A_d^{\star}(0)$. In Lemma \ref{lem:technical3}(c)-(d) we proved that $\zeta_k,\ k=0,\ldots,d-1$ are accessible fixed points in $\partial   A_d^{\star}(0)$. Consider $\zeta_0$. We know that the interval 
$(0,\zeta_0)\in \mathbb C$ determines an access to $\zeta_0$ (say, corresponding to $\beta_0$). Assume we would  have a second access to $\zeta_0$ in the upper half plane. Then by the reflection symmetry ($f(\overline{z})=\overline{f(z)}$) there would be a third access to $\zeta_0$ in the lower half plane. Using the dynamical symmetry,  this fact would reaped to each of the $\zeta_k$'s giving a total number of, at least, $3d$ accesses, a contradiction since we only have $2d$ fixed points $\beta_k$'s in $\partial D$. We conclude that $A_d^\star(0)$ has exactly $d$ accesses to infinity (and one and only one access to each of the $d$ finite fixed points in $\partial A_d^\star(0)$). Hence, $A_d^\star(0)$ is unbounded having $d$ invariant accesses to infinity.

\end{proof}

The proof of Theorem A is a direct consequence of all previous results.  We finish this section by showing that in fact the Julia set of $\Hfam$ is connected. This property is known for all Newton's method applied to any polynomial, see \cite{ConnectivityJulia}.

 \begin{theorem}
The Julia set of $\Hfam$  is a connected set in $\hat{\mathbb C}$. 
 \end{theorem}

 \begin{proof}

We need to see that all Fatou components are simply connected. We first prove this for all components corresponding to the basins of attraction of the roots. From Proposition \ref{prop:simple_connected_Halley} we need to see that eventual preimages of the immediate basins of attraction of $\{0,\alpha_1,\ldots , \alpha_{d-1}\}$ are simply connected sets in $\mathbb C$. 
From Corollary \ref{corollary:HR_coro}, all eventual preimages of the immediate attracting basins are simply connected unless one of them has more than one (free) critical point. If that was the case, using symmetry we would have that all (free) critical points should be in the same Fatou domain $U$ surrounding the origin. However, since the basins of attraction of the roots $\alpha_k$, $k=0,\ldots,d-1,$ are unbounded (see Proposition~\ref{prop:halley_alphak_unbounded}), $A_d^\star(0)$ is simply connected, and $\zeta_k\in \partial A_d^\star(0) \cap \partial A_d^\star(\alpha_k)$, $k=0,\ldots,d-1,$ then no Fatou component (other than $A_d^\star(0)$) can surround the origin. Therefore, $U$ cannot exist.

Next we see that every Fatou component not corresponding to a basin of attraction of a root is simply connected. If such a component exists, it has to be related to a critical point $c_k$, which lies in a forward invariant line $\tilde{r}_{\ell}$. This rules out the existence of rotation domains since the orbit of no critical point could accumulate on their boundaries. Therefore, such a component would correspond to a connected component of an attracting or parabolic basin of a periodic cycle contained in $\tilde{r}_{\ell}$. Moreover, the line  $\tilde{r}_{\ell}$ can contain at most two critical points $c_{k}$ and $c_{k^{\prime}}$ where $k^{\prime}=k+d/2 \ \mbox{\rm mod} (d)$, which are separated by $0$ on $\tilde{r}_{\ell}$. It is known that if a component of the basin of attraction or parabolic of an attracting or parabolic  cycle is not simply connected, then there exists a Fatou component $U$ in the basin which contains two different critical points. In particular, $U$ would contain $c_{k}$ and $c_{k^{\prime}}$. However, by the symmetry with respect to $\tilde{r}_{\ell}$, see Lemma~\ref{lem:technical2}(c), the component $U$ would surround the origin. We get to a contradiction as before as no Fatou component can surround the origin. 
 \end{proof}

\begin{figure}[htb]
	\centering
	\subfigure[\scriptsize{$d=3$ }]{\begin{tikzpicture}
			\begin{axis}[width=0.4\textwidth, axis equal image, scale only axis,  enlargelimits=false, axis on top]
				\addplot graphics[xmin=-3,xmax=3,ymin=-3,ymax=3] {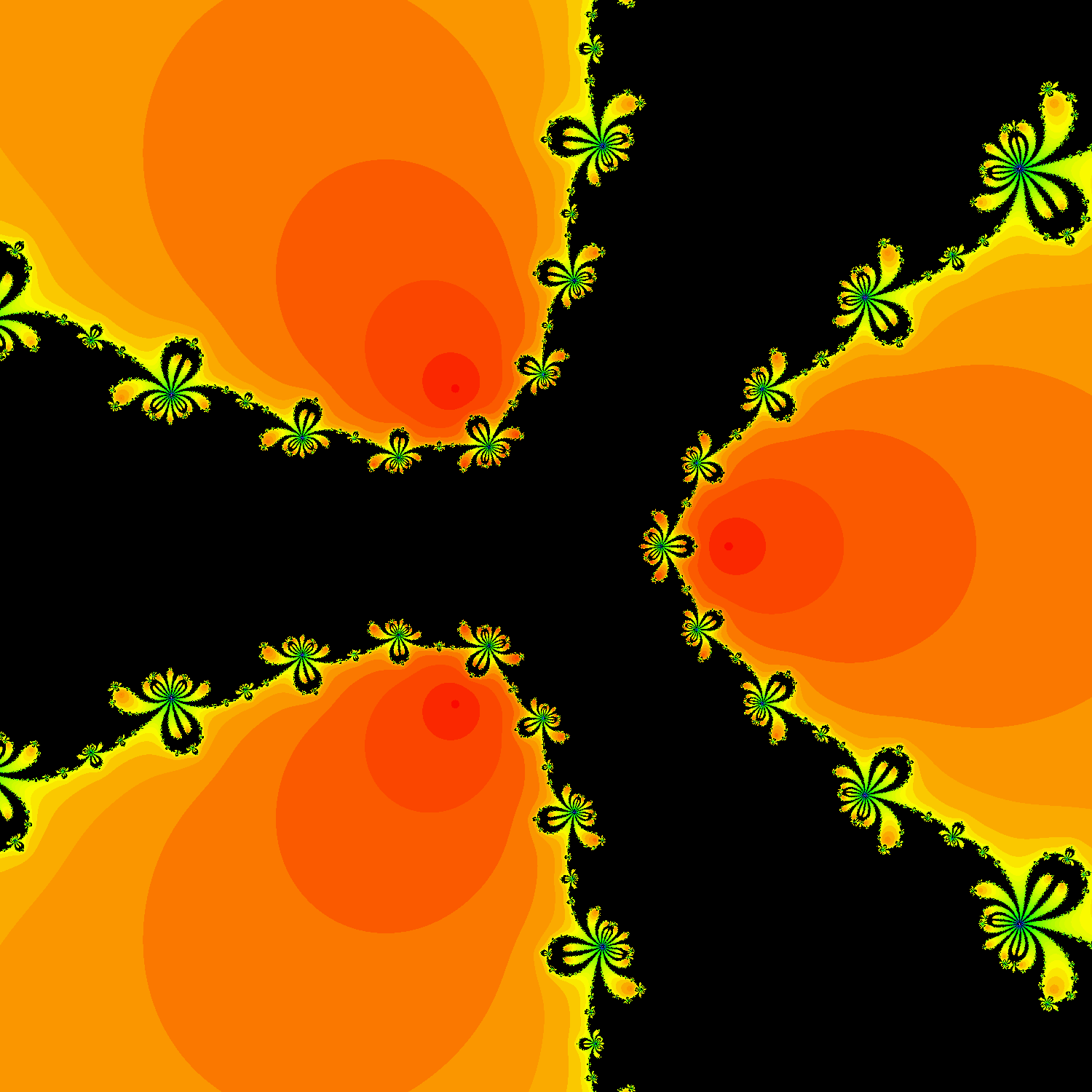};
			\end{axis}
	\end{tikzpicture}}
	\subfigure[\scriptsize{$d=4$ }  ]{	\begin{tikzpicture}
			\begin{axis}[width=0.4\textwidth, axis equal image, scale only axis,  enlargelimits=false, axis on top]
				\addplot graphics[xmin=-3,xmax=3,ymin=-3,ymax=3] {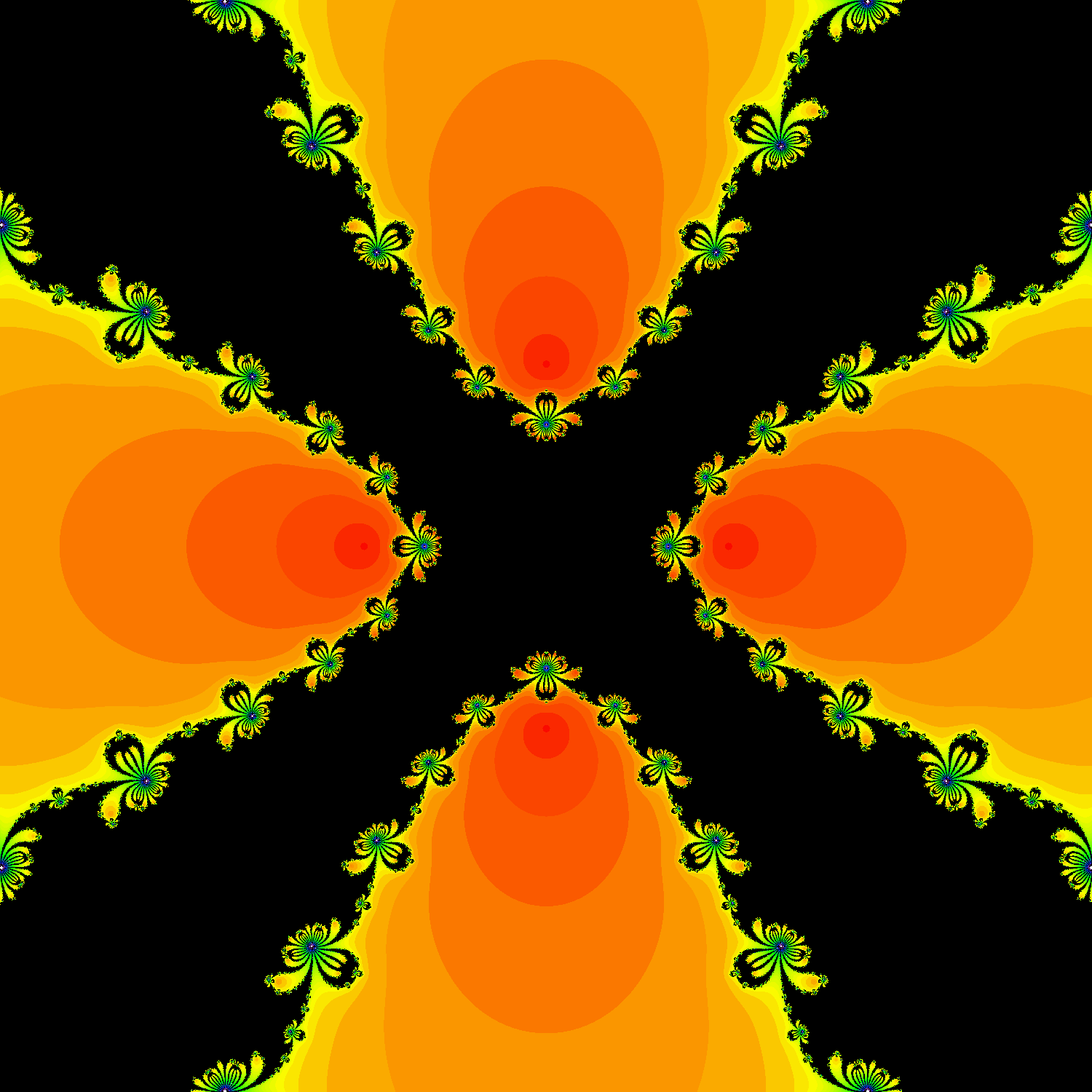};
			\end{axis}
	\end{tikzpicture}}
	\subfigure[\scriptsize{$d=5$}  ]{	\begin{tikzpicture}
			\begin{axis}[width=0.4\textwidth, axis equal image, scale only axis,  enlargelimits=false, axis on top]
				\addplot graphics[xmin=-3,xmax=3,ymin=-3,ymax=3] {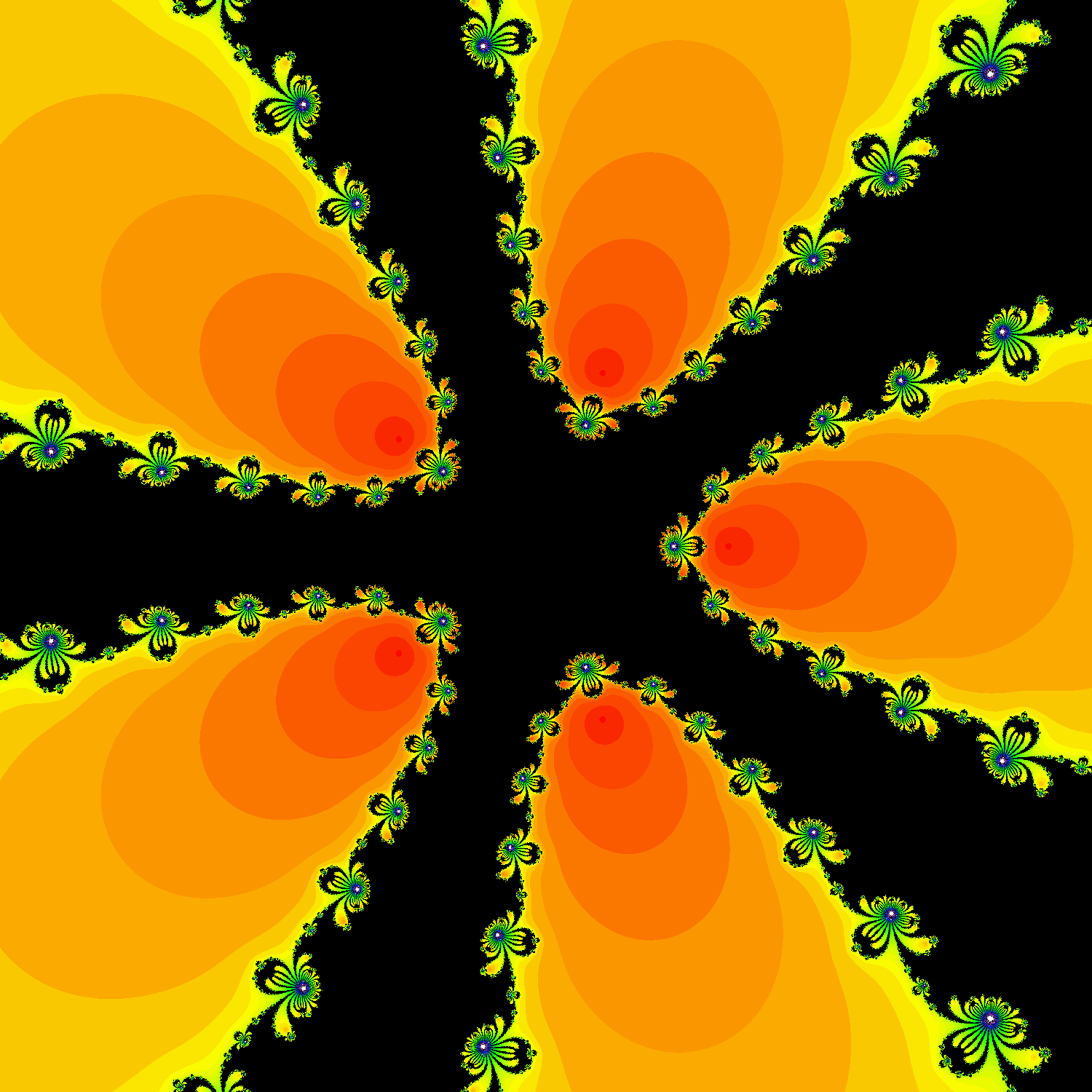};
			\end{axis}
	\end{tikzpicture}}
	\subfigure[\scriptsize{$d=6$}  ]{ \begin{tikzpicture}
			\begin{axis}[width=0.4\textwidth, axis equal image, scale only axis,  enlargelimits=false, axis on top]
				\addplot graphics[xmin=-3,xmax=3,ymin=-3,ymax=3] {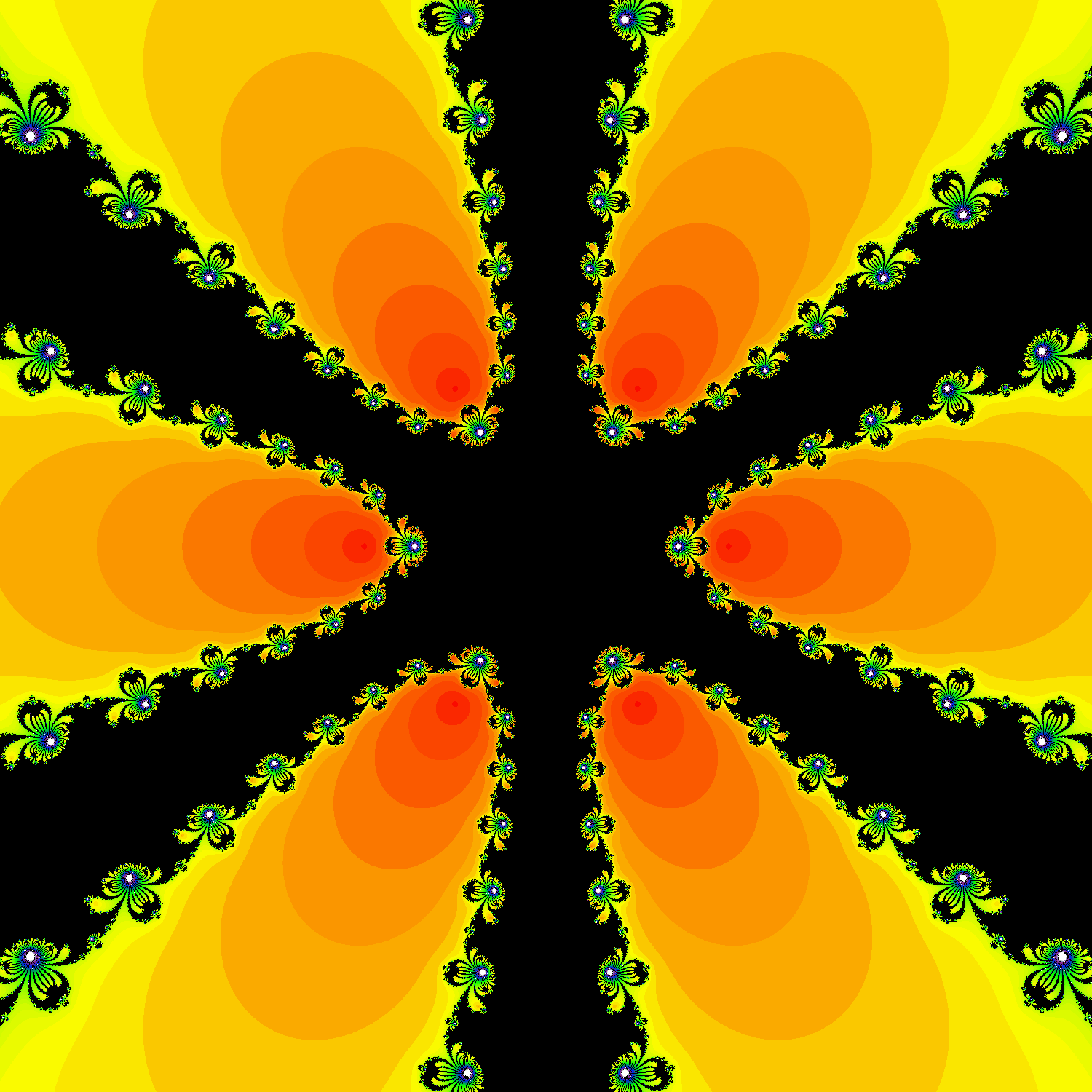};
			\end{axis}
	\end{tikzpicture}}
	\caption{\small{Dynamical planes of Traub's method applied to the family of polynomials $p_d=z(z^d-1)$. The scaling of colours is used to plot the basins of attraction of the $(d)$th-roots of the unity.}}
	\label{fig:TraubD}
\end{figure}

\section{Proof of Theorem B} \label{section:TheoremB}

Let $p$ be a polynomial of degree $d$ with all roots being simple.  Traub's method is given by the iteration $z_{n+1}=T_p(z_n)$. Traub's map (see \cite{Traub_Book}) applied to $p$ is given by the formula
\begin{equation}\label{eq:traub}
T_p (z)  =  N_p(z)  - \frac{ p(N_p(z))}{p'(z)}  
\end{equation}
where $N_p(z)$  is the Newton's map for the polynomial $p$. In the next lemma we summarize basic results about  the map $T_d$ obtained by applying Traub's method  to the polynomials $p_{d}(z)=z(z^{d}-1), \ d\geq2$.

As we did in the previous section, to simplify notation, we denote by $A_d(0):=A_{T_d}(0)$ and $A_d(\alpha_k):=A_{T_d}(\alpha_k)$, $k=0, \ldots, d-1$, the basin of attraction of each of the roots $p$ and  by $A_d^\star\left(0\right):=A^{\star}_{T_d}(\alpha_k)$ and $A_d^\star\left(\alpha_k\right):=A^{\star}_{T_d}(\alpha_k),\ k=0,\ldots d-1$, the corresponding immediate basins of attraction. We start by proving some (rather immediate) lemmas we will use later to prove Theorem B.

\begin{lemma} 
\label{lem:technical1_traub}
Let $p_{d}(z)=z(z^{d}-1), \ d\geq2$. Let $w\in \mathbb C$ and consider 
\begin{align*}
	& g_d(w)=\frac{d (d+1) w^2 ((d+1) w-1)^d -d^{d+1} w^{d+2}}{((d+1) w-1)^{d+2}}, \\
	& G_d(w)=d^{d+1} w^d-d^d (w-1) w^d+d (w-2) ((d+1)w-1)^d+(w-1) ((d+1)w-1)^d. 
\end{align*}
Then, the following statements about $\Tfam$ hold.
\begin{enumerate} [label=(\mbox{\alph*})]
\item The map $\Tfam$ is given by
$$
\Tfam(z)=\frac{d(d+1)z^{2d+1}[(d+1)z^d-1]^{d}-d^{d+1}z^{(d+1)^2}}{[(d+1)z^d-1]^{d+2}}=z g_d(z^{d}).
$$
\item The derivative of $\Tfam$ is given by
$$
\Tfam^{\prime}(z)=g_d(z^{d})+dz^{d} g_d^{\prime}(z^{d})= \frac{d (d+1) z^{2d}}{((d+1) z^d-1)^{d+3}}G_d(z^d).
$$
\item The point $z=0$ is a fixed point. Moreover, since it is also a critical point of multiplicity $2d$, it is a  super-attracting fixed point with local degree $2d+1$. 
\item The $d$th-roots of the unity are  super-attracting fixed points with local degree at least $3$ and, therefore, they are critical points of multiplicity at least $2$. 
\item The poles of $T_{p_d}$, which are the solutions of $z^d=1/(1+d)$, are mapped with degree $d+2$ onto $\infty$ and, hence, they are critical points with multiplicity $d+1$. 
\item The dynamical plane for $\Tfam$ is symmetric with respect to rotation by a $d$th-root of unity. That is, if $\eta\in \mathbb C$ such that $\eta^{d}-1=0$ and $\phi(z)=\eta z$, then
$$
T_{p_d}(z)=\left(\phi\circ T_{p_d} \circ \phi^{-1}\right)(z).
$$
\item Consider the semi-lines $r_\ell:=\{z=t\exp(\ell\pi i / d), \ t>0,\  \ell=0,\ldots, 2d-1\}$. If $\ell=2k$ for some $k\geq 0$ then  $\Tfam(r_\ell) \subset r_\ell \cup r_{\ell^{\prime}}\cup\{0\}$ with $\ell^{\prime}=(\ell + d) \ \mbox{\rm mod} (2d)$. If $\ell=2k+1$ for some $k\geq0$ then $\Tfam(r_\ell) =  r_\ell$.
\end{enumerate}
\end{lemma}

\begin{proof}
 The first equality in statement (a) follows from \eqref{eq:traub} substituting $p(z)$ by the family of polynomials $p_d(z)=z(z^d-1)$. The second equality in (a) and statement (b) follow from some computations. Statement (c) is direct from statement (b) and the fact that $G_d(0) \ne 0$. Statement (d) follows from general properties of Traub's method: fixed points corresponding to (simples) zeros of the polynomial are all critical points of degree at least 2. Alternatively one can check that $G_d(1)=G_d^{\prime}(1)=0$ and apply the dynamical plane symmetry (see statement (f)). Statement (e) follows from statement (a). Statement (f) (see statement (b) in Lemma \ref{lem:technical1}) follows from the fact that if $\eta^{d}-1=0$ then $\eta\Tfam(z)=\Tfam(\eta z)$, a direct computation from statement (a). 
 
To prove statement (g) we notice that from statement (a) we have 
$$
\Tfam(t\exp(\ell\pi i / d))=\exp(\ell\pi i / d) \ [t g_d(t^d \exp(\ell\pi i))]=\exp(\ell\pi i / d) \ [t g_d(t^d (-1)^{\ell})].
$$

Assume first that $\ell$ is even, that is $\ell=2k$ for some $k\geq 0$. Then
$$
\Tfam(t\exp(2k\pi i / d))=\exp(2k\pi i / d) \ [t g_d(t^d)].
$$ 
Clearly the function $\{t \to tg_d(t^d), \ t>0\}$ is real with one (simple, positive) pole at 
$$
t^{\star} = \left(\frac{1}{d+1}\right)^{1/d}>0.
$$
Hence $\Tfam(r_\ell) \subset r_\ell \cup r_{\ell^{\prime}}\cup \{0\}$ with $\ell^{\prime}=(\ell + d) \ \mbox{\rm mod} (2d)$.

Secondly, assume that $\ell$ is odd, that is $\ell=2k+1$ for some $k\geq 0$. Then
\begin{equation*}
\begin{split}
\Tfam(t\exp((2k+1)\pi i / d))&=\exp((2k+1)\pi i / d) \ [t g_d(-t^d)]\\
&=\exp((2k+1)\pi i / d) t \frac{d (d+1) t^{2d} ((d+1) t^d+1)^d -d^{d+1} t^{d^2+2d}}{((d+1) t^d+1)^{d+2}} \\
&:=\exp((2k+1)\pi i / d) t R_d(t).
\end{split}
\end{equation*}
Since $d(d+1)(d+1)^d-d^{d+1}>0$ we have that all coefficients of $R_d$ (numerator and denominator) are positive and hence $R_d(t)>0$ for $t>0$. We claim that $\Tfam(r_\ell) =  r_\ell$. Indeed, $\Tfam(0) = 0$ and 
\begin{equation}\label{eq:r_ell_infinity}
\lim\limits_{t \to +\infty} t R_d(t) = \lim\limits_{t \to +\infty} a_d t = +\infty \quad \mbox{with} \quad a_d=\frac{d[(d+1)^{d+1}-d^d]}{(d+1)^{(d+2)}}, \ \mbox{\rm where} \ 0<a_d<1.
\end{equation}

\end{proof}

The next proposition was proved in \cite{DavidRosado}. Nevertheless we add here a simplified proof in sake of completeness. 

\begin{proposition}
\label{prop:0-unbounded}
For $d\geq 2$ we have that $A_d^{\star}(0)$ is unbounded and it has $d$ accesses to infinity.
\end{proposition}

\begin{proof}
To prove the proposition it is enough to see that if $\ell \leq 2d-1$ is odd then $r_\ell \subset A_d^\star(0)$. From Lemma \ref{lem:technical1_traub}(g) we know that $\Tfam(r_\ell) =r_\ell$. We also know that $\Tfam(0) =0$.  Using the notation in the proof of Lemma~ \ref{lem:technical1_traub}, we have  $\Tfam(t\exp((2k+1)\pi i / d))=\exp((2k+1)\pi i / d) t R_d(t)$.  We claim that $0<tR_d(t)<t$. The first inequality is direct from the arguments used at the end of the previous proposition. The second inequality is equivalent to show that 
$R_d(t)<1$. Indeed, writing $w=t^d$  the inequality  writes as 
\begin{equation}\label{eq:T_menor_z}
P_d(w):=d (d+1) w^2 ((d+1) w+1)^d -d^{d+1} w^{d+2}-((d+1)w+1)^{d+2}<0.
\end{equation}
But this is true since a direct computation shows that all coefficients of $P_d$ are negative.

Consequently, if we denote by $S_d(t):=tR_d(t) $,
then for all $t_0>0$ the sequence $\{t_n=S_d^n(t_0)\}$ is bounded below by zero and strictly decreasing. Hence $t_n\to 0$ and $r_\ell \subset A_d^\star(0)$. 
\end{proof}

\begin{proposition}\label{prop:roots-unbounded}
 The immediate basins of attraction $A_d^{\star}(\alpha_j),\ j=0,1,\cdots, d-1$, are unbounded. Furthermore, the map $\Tfam$ has no real critical point greater than 1 and  $(1,+\infty)\subset A_d^{\star}(1) $.
\end{proposition}

\proof 
By the symmetry with respect to $d$th-roots of the unity, in order to prove that the immediate basins of attraction $A_d^{\star}(\alpha_j),\ j=0,1,\cdots, d-1$, are unbounded it is enough to prove that $A_d^{\star}(1)$ is unbounded. We will prove that $(1,+\infty)\subset A_d^{\star}(1) $. To do so it is enough to see that $\Tfam(1)=1$ and, for all $z>1$, we have $\Tfam^{\prime}(z)>0$ and $\Tfam(z)< z$. In particular, we will prove that $\Tfam$ has no real critical point greater than 1.

It is immediate that $\Tfam(1)=1$. The proof of $\Tfam(z)< z$ is equivalent to $g_d(w)<1$ which, in turn, is equivalent to \eqref{eq:T_menor_z}. Hence, to finish the argument we need to see that $T'_{p_{d}}(z)>0$ for all real $z>1$. Consider the expression of $\Tfam^{\prime}$ given in Lemma~\ref{lem:technical1_traub}(b). 
Observe that 
$$
\frac{d (d+1) z^{2d}}{((1+d) z^d-1)^{d+3}}>0 \quad \mbox{ for all } z>1.
$$
 Therefore, checking that $T'_{p_{d}}(x)>0$ for  $x>1$ is equivalent to see that the degree $d+1$ real polynomial
\begin{equation}\label{eq:g_mathcal}
\mathcal{G}_d(x)=G_d(x+1)=d^d (d-x) (x+1)^d+(d (x-1)+x) ((d+1)x+d)^d
\end{equation}
is strictly positive for $x>0$.  Easy computations show that 
\[
\mathcal{G}^{\prime}_d(x)=d^{d} (x+1)^{d-1} ( d^2-(d+1)x-1)+ d(d+1)(d(x-1)+x)((d+1)x+d)^{d-1}+ (d+1) ((d+1)x+d)^d.
\]

Notice that $x=0$ is a double zero of this polynomial (the fixed points of Traub's method corresponding to zeros of the polynomial are critical points of multiplicity at least two), but it is not easy to factor it out. Indeed, it is easy to verify that $\mathcal{G}_d(0)=\mathcal{G}^{\prime}_d(0)=0$.

We will show that the $n$th derivatives ($2\leq n \leq d+1$) of $\mathcal{G}_d$ at $x=0$ are positive, and so 
$$
\mathcal{G}_d(x) = \sum_{n=2}^{d+1}\frac{ \mathcal  G_d^{n)}(0)}{n!}   x^n>0 \ \mbox{if} \  x>0.
$$

In order to compute the $n$th derivative of $\mathcal{G}_d(x)$ we will use (a special case of) the Leibnitz formula; that is, if $u$ and $v$ are smooth functions then  
$$ 
( u v ) ^{n)} = \displaystyle \sum_{k=0}^n   
\left( 
\begin{array}{c}
n \\
k
\end{array}
 \right) 
u ^{k)} v^{n-k)}.
$$ 
In the special case where  $u$ is a polynomial of degree 1, we have that 
\begin{equation}\label{eq:Leibnitz_special}
(uv)^{n)} =  uv^{n)} + n u' v^{n-1)}.
\end{equation}
Observe that \eqref{eq:g_mathcal} can be written as
$$
\mathcal G_d(x)=u_1(x)v_1(x)+u_2(x)v_2(x) 
$$
where $u_1(x)=d^d (d-x)$, $v_1(x)=(x+1)^d$, $u_2(x)=x(d+1)-1$ and  $v_2(x)=((d+1)x+d)^d$. Hence, applying \eqref{eq:Leibnitz_special} and doing some computations we have that for $n\geq 2$
$$
 \mathcal{G}_d^{n)}(x)=\prod_{k=0}^{n-2}(d-k)\left[A_{d,n}(x)+B_{d,n}(x)\right]
$$
where
\begin{equation*}
\begin{split}
&A_{d,n}(x)=d^d(x+1)^{d-n}\left[(d-n+1)(d-x)-n(x+1)\right], \\
&B_{d,n}(x)=(d+1)^n((d+1)x+d)^{d-n}\left[(d-n+1)((d+1)x-d)+n((d+1)x+d)\right].
\end{split}
\end{equation*}

Remember that we already know that $\mathcal{G}_d^{n)}(0)=0$ for $n=0,1$.  Substituting at $x=0$ and doing some direct computations we have that  for $n\geq 2$
{\small
\begin{equation*}
\begin{split}
\mathcal{G}_d^{n)}(0)&= \prod_{k=0}^{n-2}(d-k)\left[A_{d,n}(0)+B_{d,n}(0)\right] = \\
&=d^{d-n+1}\prod_{k=0}^{n-2}(d-k) \left[d^{n+1}+(1-n)d^n-nd^{n-1}+\sum_{j=0}^{n} (2n-1)\left(\begin{array}{c}n \\ j\end{array}\right)d^j-\sum_{j=0}^{n}\left(\begin{array}{c}n \\ j\end{array}\right)d^{j+1}\right]  \\
&:=d^{d-n+1}\prod_{k=0}^{n-2}(d-k) \ C_{d,n}(0).
\end{split}
\end{equation*}
}
Notice that $k\leq n-2 \leq d-1<d$. So, to show that $\mathcal{G}_d^{n)}(0)\geq 0$ for $n\geq 0$ it is enough to show that $C_{d,n}(0)\geq 0$ for $n\geq 2$. One can easily check from its definition that, for all $n=2,\ldots, d+1$, the expression of $C_{d,n}(0)$ is a degree $n+1$ polynomial in the variable $d$, that is,
$$
C_{d,n}(0)=\sum_{\ell=1}^{n+1} c_{\ell}(d,n) d^{\ell}, \quad d\geq 2, \ n\leq d+1. 
$$
We claim that $c_\ell(d,n)\geq 0$ concluding the desired result $\mathcal{G}_d^{n)}(0)$ for all suitable values of $d$ and $n$. Indeed, given $d\geq 2$ and $n\leq d+1$ we have $n+2$ coefficients in the variable $d$. Moreover,  
\begin{equation*}
\begin{split}
&c_{n+1}(d,n)=c_{n}(d,n)=0,\ c_{n-1}(d,n)=\frac{3}{2}n(n-1)>0,\ c_0(d,n)=2n-1>0, \mbox{ and }\\
&c_\ell(d,n)=(2n-1)\frac{n!}{(n-\ell)!\ell!}-\frac{n!}{(n-(\ell-1))!(\ell-1)!}>0 \quad \mbox{ for } 1\leq \ell \leq n-2. 
\end{split}
\end{equation*}

All together implies that $[1,\infty) \subset A_d^{\star}(1)$ and so $A_d^{\star}(1)$ is unbounded. By the dynamical symmetry  described in Lemma \ref{lem:technical1_traub} the statement of the proposition is proved.

\endproof

The next Lemma provides a criterion for the connectivity of Julia sets. It is particularly interesting because it does not make use of  fixed points of the Riemann-Hurwitz formula. Instead, it uses the existence of a collection of unbounded Fatou components that contain in their boundaries all poles.  Similar arguments can be found in the literature (see, for instance, \cite[Lemma 4.5]{GNPP}).

\begin{lemma}\label{lem:unbounded}
Let $R$ be a Rational map and let $\mathcal{U}=\{U_1,...,U_{\ell}\}, \ \ell\geq2$, a  collection of unbounded Fatou components of $R$.  Assume that if $p$ is a pole of $R$ then $p\in \partial U_{i}\cap \partial U_{j}$ for some $i,j\in\{1,\ldots,\ell\}$ with $i\ne j$. Then, the Julia set of $R$ is connected.
\end{lemma}

\begin{proof}
First observe that since $\infty$ belongs to the boundary of $\ell\geq2$ Fatou components, then $\infty$ belongs to the Julia set.
 
Assume that the Julia set of $R$ is not connected. Then, there is a multiply connected Fatou component $V$ and we can take a Jordan curve $\gamma\subset V$ such that $\gamma$ separates a (bounded) component $\Gamma$ of the Julia set  and $\infty$. Since backwards preimages of $\infty$ are dense in the Julia set, then the component $\Gamma$ must contain pre-poles. Thus, there is a minimal $m\geq 0$ such that $V_m:=R^m(V)\subset \mathcal F(R)$ separates a pole $p$ from $\infty$. However, this is impossible since, by assumption, $p$ belongs to the boundary of two different  (unbounded) Fatou components in $\mathcal{U}$ (so, one of them cannot be $V_m$).  We have reached a contradiction and, therefore, we can conclude that the Julia set of $R$ is connected.

\end{proof}

The next theorem is a consequence of Lemma~\ref{lem:unbounded} and, in particular, implies Theorem B. 

\begin{theorem}\label{thm:connectivityTraub}
	The Julia set of $\Tfam$ is connected.
\end{theorem}
\begin{proof}
	Let $\mathcal{U}=\{A_d^{\star}(0),A_d^{\star}(\alpha_1),\ldots  A_d^{\star}(\alpha_d)\}$. By Proposition~\ref{prop:0-unbounded} and Proposition~\ref{prop:roots-unbounded}, $\mathcal{U}$ is a collection of unbounded Fatou components. By Lemma~\ref{lem:unbounded}, to finish the proof it is enough to show that every pole of $\Tfam$ belongs to the boundary of, at least, two different Fatou components in $\mathcal U$.
	
The poles of $\Tfam$ are the solutions of $z^d=1(1+d)$. Therefore, $\Tfam$ has exactly $d$ poles, denoted by $p_j,\ j=1,\ldots, d$. We claim that, after renumbering if necessary, we have  
\begin{equation}\label{eq:poles}
p_j \in \partial A_d^{\star}(0) \cap A_d^{\star}(\alpha_j), \ j=1,\ldots, d.
\end{equation}
	
We first prove that each pole belongs to the boundary of the immediate basin of attraction of a $d$th-root of the unity (we use this fact to denote each pole by $p_j$, accordingly). From  Lemma~\ref{lem:technical1_traub}, it is enough to prove that there is a pole, which will be denoted by $p_1$, in $\partial A_d^{\star}(1)$. By Proposition~\ref{prop:roots-unbounded} we know that $[1,+\infty)\subset A_d^{\star}(1)$. The fixed point $z=1$ has local degree 3 and, so, the segment $(1,+\infty)$ has three preimages, say $L_k,\ k=1,2,3$, in $A_d^{\star}(1)$. One of them, which we denote by $L_1$, is itself. That is $L_1=(1,+\infty)$. Since $L_1$ lands at $z=\infty$, which is a repelling fixed point  (and therefore has local degree 1), the non-zero landing points of $L_2$ and $L_3$ cannot be  $z=\infty$. Hence, both $L_2$ and $L_3$ should be arcs in $A_d^{\star}(1)$ connecting $z=1$ with a pole. The dynamical symmetry described in Lemma~\ref{lem:technical1_traub} implies $p_j\in \partial A_d^{\star}(\alpha_j)$.
	
The fact that all poles belong to $\partial A_d^{\star}(0)$ is proved analogously using the invariant lines $r_\ell$ with $\ell$ odd in $A_d^{\star}(0)$ found in Proposition \ref{prop:0-unbounded}.  Assume $\ell=1$ (the other lines follow by the symmetry). Since the local degree of $\Tfam$ at the origin is $2d+1$ there must be at least $2d+1$ preimages of $r_1$ in $A_d^{\star}(0)$ - say $L_k,\ k=1,\ldots 2d+1$ but only one might land at infinity (since $z=\infty$ has local degree $1$). Hence, by symmetry, there must be at least one arc connecting $z=0$ and each pole $p_j, \ j=1,\ldots, d$.
\end{proof}

\bibliographystyle{alpha}
\bibliography{biblio}

\end{document}